\documentclass{article}

\usepackage{float}
\usepackage{xcolor}
\usepackage{arxiv}
\usepackage{amsmath}
\usepackage[utf8]{inputenc} 
\usepackage{textcomp}

\usepackage[T1]{fontenc}    
\usepackage{hyperref}       
\usepackage{url}            
\usepackage{booktabs}       
\usepackage{amsfonts}       
\usepackage{nicefrac}       
\usepackage{microtype}      
\usepackage{lipsum}
\usepackage{tikz}
\usepackage{stmaryrd}
\usepackage{amssymb}
\usepackage{amsthm}
\usepackage{mathtools}
\usetikzlibrary{calc, positioning}
\usepackage{esdiff}
\usepackage{pgfplots}
\pgfplotsset{compat=1.14}

\usetikzlibrary{backgrounds,pgfplots.groupplots}
\definecolor{ao(english)}{rgb}{0.0, 0.5, 0.0}
\usepackage{letltxmacro}
\LetLtxMacro{\originaleqref}{\eqref}
\renewcommand{\eqref}{Eq.~\originaleqref}
\usepackage{physics}
\usetikzlibrary{arrows}      
\usetikzlibrary{calc} 
\usetikzlibrary{calc,patterns,angles,quotes}
\usepackage{siunitx}

\usepackage{lipsum}
\usepackage{titlesec}
\usepackage{array}
\newcolumntype{P}[1]{>{\centering\arraybackslash}p{#1}}
\usepackage{tabularx, multirow, booktabs}
\usepackage{cite}

\usepackage[title]{appendix}

\makeatletter
\newlength{\flex@pattern@density}
\newlength{\flex@pattern@linewidth}
\newlength{\flex@pattern@auxlength}
\newlength{\flex@pattern@auxlengthtwo}
\tikzset{/tikz/.cd,
    pattern density/.code={\setlength\flex@pattern@density{#1}
    \pgfmathsetmacro{\tmp}{1.1*#1}
    \setlength\flex@pattern@auxlength{\tmp pt}
    \setlength\flex@pattern@auxlengthtwo{#1}
    \advance\flex@pattern@auxlengthtwo by 0.1pt
    \typeout{\the\flex@pattern@density,\the\flex@pattern@auxlength}},
    pattern density=3pt,
    pattern line width/.code={\setlength\flex@pattern@linewidth{#1}},
    pattern line width=0.4pt,
}

\pgfdeclarepatternformonly[\flex@pattern@density,\flex@pattern@linewidth]{flexible horizontal lines}{\pgfpointorigin}{\pgfqpoint{100pt}{1pt}}{\pgfqpoint{100pt}{\flex@pattern@density}}%
{
  \pgfsetlinewidth{\flex@pattern@linewidth}
  \pgfsetcolor{\tikz@pattern@color}
  \pgfpathmoveto{\pgfqpoint{0pt}{0.5pt}}
  \pgfpathlineto{\pgfqpoint{100pt}{0.5pt}}
  \pgfusepath{stroke}
}

\pgfdeclarepatternformonly[\flex@pattern@density,\flex@pattern@linewidth]{flexible vertical lines}{\pgfpointorigin}{\pgfqpoint{1pt}{100pt}}{\pgfqpoint{\flex@pattern@density}{100pt}}%
{
  \pgfsetlinewidth{\flex@pattern@linewidth}
  \pgfsetcolor{\tikz@pattern@color}
  \pgfpathmoveto{\pgfqpoint{0.5pt}{0pt}}
  \pgfpathlineto{\pgfqpoint{0.5pt}{100pt}}
  \pgfusepath{stroke}
}

\pgfdeclarepatternformonly[\flex@pattern@auxlengthtwo,\flex@pattern@density,\flex@pattern@linewidth,\flex@pattern@auxlength]{flexible north east lines}{\pgfqpoint{-1pt}{-1pt}}{\pgfqpoint{\flex@pattern@auxlength}{\flex@pattern@auxlength}}{\pgfqpoint{\flex@pattern@density}{\flex@pattern@density}}%
{
  \pgfsetlinewidth{\flex@pattern@linewidth}
  \pgfsetcolor{\tikz@pattern@color}
  \pgfpathmoveto{\pgfqpoint{0pt}{0pt}}
  \pgfpathlineto{\pgfqpoint{\flex@pattern@auxlengthtwo}{\flex@pattern@auxlengthtwo}}
  \pgfusepath{stroke}
}

\pgfdeclarepatternformonly[\flex@pattern@auxlengthtwo,\flex@pattern@density,\flex@pattern@linewidth,\flex@pattern@auxlength]{flexible north west lines}{\pgfqpoint{-1pt}{-1pt}}{\pgfqpoint{\flex@pattern@auxlength}{\flex@pattern@auxlength}}{\pgfqpoint{\flex@pattern@density}{\flex@pattern@density}}%
{
  \pgfsetlinewidth{\flex@pattern@linewidth}
  \pgfsetcolor{\tikz@pattern@color}
  \pgfpathmoveto{\pgfqpoint{\flex@pattern@auxlengthtwo}{0pt}}
  \pgfpathlineto{\pgfqpoint{0pt}{\flex@pattern@auxlengthtwo}}
  \pgfusepath{stroke}
}

\usepackage{setspace}

\titlespacing\section{0pt}{2pt plus 0pt minus 0pt}{2pt plus 0pt minus 0pt}
\titlespacing\subsection{0pt}{2pt plus 0pt minus 0pt}{2pt plus 0pt minus 0pt}
\titlespacing\subsubsection{0pt}{2pt plus 0pt minus 0pt}{2pt plus 0pt minus 0pt}
\pgfdeclarepatternformonly[\flex@pattern@density,\flex@pattern@linewidth,\flex@pattern@auxlength,\tikz@pattern@color]{flexible grid}{\pgfqpoint{-1pt}{-1pt}}{\pgfqpoint{\flex@pattern@auxlength}{\flex@pattern@auxlength}}{\pgfqpoint{\flex@pattern@density}{\flex@pattern@density}}%
{
  \pgfsetlinewidth{\flex@pattern@linewidth}
  \pgfsetcolor{\tikz@pattern@color}
  \pgfpathmoveto{\pgfqpoint{0pt}{0pt}}
  \pgfpathlineto{\pgfqpoint{0pt}{\flex@pattern@density}}
  \pgfpathmoveto{\pgfqpoint{0pt}{0pt}}
  \pgfpathlineto{\pgfqpoint{\flex@pattern@density}{0pt}}
  \pgfusepath{stroke}
}

\pgfdeclarepatternformonly[\flex@pattern@density,\flex@pattern@linewidth,\flex@pattern@auxlength,\tikz@pattern@color]{flexible crosshatch}{\pgfqpoint{-1pt}{-1pt}}{\pgfqpoint{\flex@pattern@auxlength}{\flex@pattern@auxlength}}{\pgfqpoint{3pt}{3pt}}%
{
  \pgfsetlinewidth{\flex@pattern@linewidth}
  \pgfsetcolor{\tikz@pattern@color}
  \pgfpathmoveto{\pgfqpoint{\flex@pattern@density}{0pt}}
  \pgfpathlineto{\pgfqpoint{0pt}{\flex@pattern@density}}
  \pgfpathmoveto{\pgfqpoint{0pt}{0pt}}
  \pgfpathlineto{\pgfqpoint{\flex@pattern@density}{\flex@pattern@density}}
  \pgfusepath{stroke}
}

\pgfdeclarepatternformonly[\flex@pattern@density,\flex@pattern@linewidth,\tikz@pattern@color]{flexible dots}{\pgfqpoint{-1pt}{-1pt}}{\pgfqpoint{1pt}{1pt}}{\pgfqpoint{\flex@pattern@density}{\flex@pattern@density}}%
{
  \pgfsetcolor{\tikz@pattern@color}
  \pgfpathcircle{\pgfqpoint{0pt}{0pt}}{\flex@pattern@linewidth}
  \pgfusepath{fill}
}

\pgfdeclarepatternformonly[\flex@pattern@density,\flex@pattern@linewidth,\tikz@pattern@color]{flexible crosshatch dots}{\pgfqpoint{-1pt}{-1pt}}{\pgfqpoint{2.5pt}{2.5pt}}{\pgfqpoint{\flex@pattern@density}{\flex@pattern@density}}%
{
  \pgfsetcolor{\tikz@pattern@color}
  \pgfpathcircle{\pgfqpoint{0pt}{0pt}}{\flex@pattern@linewidth}
  \pgfpathcircle{\pgfqpoint{1.5pt}{1.5pt}}{\flex@pattern@linewidth}
  \pgfusepath{fill}
}

\makeatother

\usetikzlibrary{quotes,arrows.meta}
\tikzset{
  annotated cuboid/.pic={
    \tikzset{%
      every edge quotes/.append style={midway, auto},
      /cuboid/.cd,
      #1
    }
    \draw [every edge/.append style={pic actions, densely dashed, opacity=.5}, pic actions]
    (0,0,0) coordinate (o) -- ++(-\cubescale*\cubex,0,0) coordinate (a) -- ++(0,-\cubescale*\cubey,0) coordinate (b) edge coordinate [pos=1] (g) ++(0,0,-\cubescale*\cubez)  -- ++(\cubescale*\cubex,0,0) coordinate (c) -- cycle
    (o) -- ++(0,0,-\cubescale*\cubez) coordinate (d) -- ++(0,-\cubescale*\cubey,0) coordinate (e) edge (g) -- (c) -- cycle
    (o) -- (a) -- ++(0,0,-\cubescale*\cubez) coordinate (f) edge (g) -- (d) -- cycle;
    \path [every edge/.append style={pic actions, |-|}];
    
  },
  /cuboid/.search also={/tikz},
  /cuboid/.cd,
  width/.store in=\cubex,
  height/.store in=\cubey,
  depth/.store in=\cubez,
  units/.store in=\cubeunits,
  scale/.store in=\cubescale,
  width=10,
  height=10,
  depth=10,
  units=cm,
  scale=.1,
}
\makeatletter
\pgfmathdeclarefunction{erf}{1}{%
  \begingroup
    \pgfmathparse{#1 > 0 ? 1 : -1}%
    \edef\sign{\pgfmathresult}%
    \pgfmathparse{abs(#1)}%
    \edef\x{\pgfmathresult}%
    \pgfmathparse{1/(1+0.3275911*\x)}%
    \edef\t{\pgfmathresult}%
    \pgfmathparse{%
      1 - (((((1.061405429*\t -1.453152027)*\t) + 1.421413741)*\t 
      -0.284496736)*\t + 0.254829592)*\t*exp(-(\x*\x))}%
    \edef\y{\pgfmathresult}%
    \pgfmathparse{(\sign)*\y}%
    \pgfmath@smuggleone\pgfmathresult%
  \endgroup
}
\makeatother

\DeclareMathOperator{\DDiv}{div}

\DeclareMathOperator{\Div}{\nabla\cdot}
\newcommand{\pjump}{\llbracket p \rrbracket}
\newcommand{\Spaceu}{\tilde H(\DDiv, \Omega)}
\DeclareMathOperator{\Subdomain}{\Omega\backslash\gamma}
\newcommand{\SpacepBVP}{H^1(\Subdomain)}
\newcommand{\perm}{\mathbf{\kappa}}
\newcommand{\dx}{d\mathbf{x}}
\newcommand{\ds}{d\mathbf{\sigma}}
\newcommand{\pD}{p}
\newcommand{\vari}{(\mathbf{x})}
\newtheorem{thm}{Theorem}[section]
\newtheorem{definition}{Definition}[section]
\newcommand{\dapprox}{\delta_\epsilon^\prime}
\newcommand{\deltaeps}{\delta_{\epsilon}}
\newcommand{\pjumpinu}{\frac{u_{n,\gamma}}{t_f}}
\newcommand{\dipole}{G_{\epsilon}}
\newcommand{\Dipole}{D_{\epsilon}}
\newenvironment{myproof}[1][\proofname]{\proof[#1]\mbox{}\\*}{\endproof}

\numberwithin{equation}{section}

\definecolor{c1}{RGB}{170,68,153}
\definecolor{c2}{RGB}{136,34,85}
\definecolor{c3}{RGB}{17,119,51}
\definecolor{c4}{RGB}{51,34,136}

\title{An efficient method for modeling flow in porous media with immersed faults}

\author{
Youguang Chen \\
  Oden Institute for Computational Engineering and Sciences\\
  The University of Texas at Austin\\
  \texttt{youguang@utexas.edu} \\
   \And
George Biros \\
  Oden Institute for Computational Engineering and Sciences\\
  The University of Texas at Austin\\
  \texttt{biros@oden.utexas.edu} \\
}

\begin{document}
\maketitle

\begin{abstract}
Modeling flow in geosystems with natural fault is a challenging problem due to low permeability of fault compared to its surrounding porous media. One way to predict the behavior of the flow while taking the effects of fault into account is to use the mixed finite element method. However, the mixed method could be time consuming due to large number of degree of freedom since both pressure and velocity are considered in the system. A new modeling method is presented in this paper. First, we introduce approximations of pressure based on the relation of pressure and velocity. We further decouple the approximated pressure from velocity so that it can be solved independently by continuous Galerkin finite element method. The new problem involves less degree of freedom than the mixed method for a given mesh . Moreover, local problem associated with a small subdomain around the fault is additionally solved to increase the accuracy of approximations around fault. Numerical experiments are conducted to examine the accuracy and efficiency of the new method. Results of three-dimensional tests show that our new method is up to 30$\times$ faster than the the mixed method at given $L^2$ pressure error.

\end{abstract}

\section{Introduction}

A fault is a fracture in a volume of rock which has lower permeability than its surrounding matrix. Since the permeability difference can vary several magnitudes, a fault has significant effects on flow in porous media by acting either as a conduit causing flow retardation or as a barrier restricting flow going through it. Predicting flow behavior for geological systems with fault zones is important in many applications, such as exploitation of oil and geothermal resources \cite{geothermal}, and CO\textsubscript{2}-sequestration \cite{antonio}. 

The influences of faults on flow are complicated since many petrophysical properties are involved, but these effects can be simplified to derive fault models that are suitable for simulations. These models for faults can be classified as continuum models and discrete models. Continuum models treat faults implicitly as the same with its surrounding rocks by homogenizing processes, while the discrete models can consider the faults individually. We restrict our attention to the discrete model in this paper because the flow transport is not adequately captured in the continuum models. Since the width of fault is very small compared to the characteristic length of the whole simulation domain, faults are often represented as $(d-1)-$dimensional immersed interfaces in $d-$dimensional domain. Under such case, refined grids inside the fault can be avoided and thus the computational costs could be saved. Such reduced model for fractures was derived by Alboin et al. \cite{alboin2002} for single-phase Darcy flow in porous media by coupling conditions at the fracture-matrix interface. Later, this model was extended by Martin et al. \cite{martin2005} for faults by considering more general coupling conditions. Moreover, \cite{hoteit2008} and \cite{jaffre2011} extended this model for applications in two-phase flow problems.

In this paper, we study the model presented in \cite{martin2005} where pressure could have jump while normal component of velocity is continuous along the fault. We present a new formulation that could be solved by the continuous Galerkin method to generate approximations of the pressure. Based on this formulation, a new method is proposed to generate pressure and velocity solutions to the boundary value problem. 

\textbf{Contributions.} We summarize our contributions as follows:
\begin{itemize}
    \item We define the boundary value problem of flow in porous media with fault and introduce approximations of pressure, which is based on the relation of Laplacian of pressure to the velocity along the fault.
    \item To decouple pressure and velocity, we derive equivalent model in one-dimensional case, and approximate models when dimension is larger than one. The new formulation enables to solve pressure using continuous Galerkin finite element method.
    \item To correct approximations of the solutions around faults, we define and solve problems associated with small subdomains around faults using the mixed finite element method. The subdomain problems could be independently solved for each fault. 
\end{itemize}

\textbf{Limitations.} One limitation is that we assume uniform normal direction for each fault (point in $\mathbb{R}^1$, line in $\mathbb{R}^2$ and plane in $\mathbb{R}^3$) in our derivations of the new model. But faults with complicated geometries could be decomposed at first to faults with uniform normal directions, and then the problems could be solved using the method proposed in this paper. Another limitation of our method is that we need to solve the flow problem using the mixed method with subdomains associated with faults in order to get accurate solutions in the whole domain. The number of degree of freedom may be large in the case of three-dimensional space with some intersected faults, and thus could be time consuming.

\textbf{Related work.}
 Both finite volume and finite element methods have been proposed as discretization methods for the discrete model of flow in porous media with immersed faults. In \cite{h2009}, the authors introduce both vertex-centered and cell-centered finite volume methods in two-dimensional domain to account for the flow between faults and matrix. In \cite{faille2016} and \cite{fumagalli2018}, the authors propose finite volume methods for the application of non-matching grids along faults. For the finite element approaches dealing with flow simulations with immersed faults, mixed finite element method is robust to handle the pressure discontinuity occurred along the fault and satisfy the mass conservation property. In \cite{martin2005}, a  mixed finite element scheme for the coupled pressure-velocity system of the reduced discrete model where faults are treated as interfaces. In \cite{carlo2012,fumagalli2011,fumagalli2013,fumagalli2014,marco2017}, the mixed finite element method was extended for non-matching grids by using extended finite element method. Using non-matching grids can simplify the mesh for realistic problems and make it possible to run multiple cases without remeshing for different fault configurations, which is useful in the study of uncertainty quantification of geophysical properties. Furthermore, mortar technique was applied to the mixed finite element method to solve flow with non-matching grids with respect to the complex networks of fractures or faults in \cite{pichot2010,pichot2012,boon2018}.

\textbf{Outline of the paper.} In \autoref{formulation} we present the boundary value problem considered in this paper with the corresponding variational problem of the mixed method. Moreover, we define the approximate pressure and derive differential equations about it. In \autoref{implementation} we remark some implementation issues related to the new method. In \autoref{results} numerical experiments are performed under different mesh resolutions for both the mixed method and the new method, and we report and compare the accuracy of the solutions and the CPU time used for each test.

\section{Formulation}\label{formulation}

\subsection{Boundary value problem}

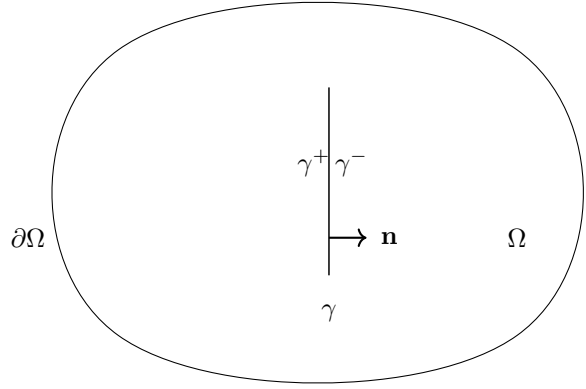
\begin{figure}[ht!]
\begin{center}
  \begin{tikzpicture}
     \path[font={\tiny}]
        (11 , -.8)   coordinate (A)
        (5.7   , -.8)   coordinate (B)
        (5.7   , 3)   coordinate (C)
        (11   , 3)   coordinate (D)
    ;
    \draw plot [smooth cycle, tension=0.8] coordinates {(A) (B) (C) (D)};
    \draw[black, line width = 0.2mm] (8.5,0) to (8.5,2.5);
            \node[] at (11,0.5)   (a) {$\Omega$};
        \node[] at (4.5,0.5)   (a) {$\partial\Omega$};
        \node[] at (8.5, -0.5)   (b) {$\gamma$};
        \node[] at (8.3, 1.5) (c) {$\gamma^+$};
        \node[] at (8.8,1.5) (d) {$\gamma^-$};
        \node[] at (9.3,0.5) (d) {$\mathbf{n}$};
    \draw[thick, ->] (8.5, 0.5) to (9,0.5);
  \end{tikzpicture}
  \end{center}
  \caption{The domain $\Omega$ with an immersed fault $\gamma$.}
  \label{domain}
\end{figure}
We attempt to consider a bounded domain $\Omega \in \mathbb{R}^d$ ($d$ = 1, 2 or 3) with Lipschitz boundary denoted by $\partial \Omega$. Within the domain, we consider a fault $\gamma \in \mathbb{R}^{d-1}$ with normal direction denoted by $\mathbf{n}$. We denote the side of $\gamma$ with outward normal direction $\mathbf{n}$ as $\gamma^+$, and the other side as $\gamma^-$. \autoref{domain} shows an example of such configuration.\par

For the flow in $\Omega$, it is supposed to satisfy the conservation equation and Darcy's law:
    \begin{align}\label{BVP_domain}
    \begin{cases}
             \Div \mathbf{u} = f \qquad &\text{in $\Omega$},  \\
 	     \mathbf{u} = -\frac{\rho}{\mu}  \mathbf{\kappa} \nabla p \qquad &\text{in $\Omega\backslash\gamma$},   
    \end{cases}
    \end{align}  
where $\mathbf{u}$ is the velocity and $p$ is the pressure, $f$ is the external source, $\mathbf{\kappa}$ is the permeability tensor of the domain, and $\rho$, $\mu$ are the density and viscosity of the fluid, respectively. We assume that $f\in L^2(\Omega)$,  $\perm \in (L^\infty(\Omega))^{(d\times d)}$ with positive components, $\rho$ and $\mu$ are positive constants. \par

For the fault, we denote its transmissibility by $t_f\in\mathbb{R^+}$, a parameter indicating the ability of transporting fluids through the fault. We remark that we treat the fault with no thickness, i.e., $\gamma \in \mathbb{R}^{d-1}$. The transmissibility can be derived by
        \begin{align*}
		t_f = \frac{\mathbf{\kappa}_f}{d_f},
        \end{align*}
where $\mathbf{\kappa}_f$ and $d_f$ are the permeability and the thickness of the fault. Since the fault has different permeability compared with its neighboring domain, it generally causes pressure jump $\pjump$ along $\gamma^+$ and $\gamma^-$ which is defined by 
        \begin{equation}\label{p_jump}
		\pjump = \gamma_0^+p - \gamma_0^-p,
        \end{equation}
where $\gamma_0^+$ and $\gamma_0^-$ are trace operators defined on the side of $\gamma^+$ and $\gamma^-$, respectively.

Moreover, we have a constitutive condition of the velocity and the pressure along the fault:
	\begin{align}\label{BVP_fault}
	 \gamma_0^+( \mathbf{u}\cdot\mathbf{n}) = \gamma_0^-(\mathbf{u}\cdot\mathbf{n})= t_f \pjump \qquad \text{on $\gamma$}.
	\end{align}  
	\par

To make sense of the previous boundary value problem, we need to define spaces where our pressure and velocity solutions should lie. From the first equation of \eqref{BVP_domain}, we naturally expect the velocity at least from $H(\DDiv, \Omega) = \{\mathbf{v}\in (L^2(\Omega))^d: \Div \mathbf{v} \in L^2(\Omega)\}$. Indicated by the second equation of \eqref{BVP_domain}, $p\in H^1(\Omega \backslash \gamma)$. Then \eqref{p_jump} is well defined in the sense that $\gamma_0^+p$, $\gamma_0^-p \in H^{1/2}(\gamma)$ by trace theorem. The first equality of \eqref{BVP_fault} indicates that normal component of the velocity along the fault should be continuous, which is not generally satisfied in $H(\DDiv, \Omega)$. Besides, in general $\gamma_0^+ \text{ and } \gamma_0^-: \mathbf{u}\cdot \mathbf{n} \xrightarrow{} H^{-1/2}(\gamma)$ for $\mathbf{u} \in H(\DDiv,\Omega)$. Thus, we need more restrictions for the velocity such that the traces of the pressure and the normal component of the velocity could be consistent indicated by the second equality of \eqref{BVP_fault}. Then we can get the solution space for velocity to
\begin{align*}
    \Spaceu = \{\mathbf{v} \in H(\DDiv,\Omega): \gamma_0^+( \mathbf{u}\cdot\mathbf{n}) = \gamma_0^-(\mathbf{u}\cdot\mathbf{n})\in H^{1/2}(\gamma) \text{ a.e. on $\gamma$} \}.
\end{align*}
It is easy to verify that the new defined space for velocity is a closed subspace of $H(\DDiv, \Omega)$, and thus a Hilbert space. We define our solution space as
\begin{align*}
    \mathcal{H} = \Spaceu \times \SpacepBVP.
\end{align*}

For simplicity of notation, we can define 
\begin{align*}
    \mathbf{u}\cdot\mathbf{n} =  \gamma_0^+( \mathbf{u}\cdot\mathbf{n}) = \gamma_0^-(\mathbf{u}\cdot\mathbf{n}) \qquad \forall \mathbf{u} \in \Spaceu.
\end{align*}

We additionally assume that $\frac{\rho}{\mu} = 1$ and Dirichlet boundary condition for pressure is imposed on $\partial \Omega$. Now the boundary value problem is well defined and we formulate it as the follows: Find $(\mathbf{u}, p) \in \mathcal{H}$ such that
	  \begin{align}\label{BVP}
          \begin{cases}
                   \Div \mathbf{u} = f \qquad &\text{in $\Omega$},  \\
                  \mathbf{\kappa}^{-1} \mathbf{u} +  \nabla p  = 0 \qquad &\text{in $\Omega\backslash\gamma$},   \\
		   t_f^{-1}\mathbf{u} \cdot \mathbf{n}  = \pjump  &\text{on $\gamma$}, \\
                   p = p_{\text{D}} &\text{on $\partial \Omega$}.
          \end{cases}
          \end{align}


\subsection{Mixed variational problem}
In order to get the mixed variational formulation, we first take dot product of the second equation of \eqref{BVP} by $\mathbf{v} \in \Spaceu$ and use integration-by-parts
	\begin{align*}
		\int_{\Omega} \perm^{-1}\mathbf{u} \cdot \mathbf{v} \, \dx \nolinebreak
		&=  -\int_{\Subdomain} \nabla p \cdot \mathbf{v} \, \dx \nolinebreak
		 = \int_{\Subdomain} p \Div \mathbf{v} \, \dx - \int_{\partial(\Subdomain)} p \mathbf{v} \cdot \nu \, \dx \\
		&= \int_{\Omega} p \Div \mathbf{v} \, \dx \nolinebreak
		    - \int_{\partial\Omega}p_D \mathbf{v} \cdot \mathbf{\nu} \, \ds \nolinebreak
		    - \int_{\gamma^+} (\gamma_0^+p) (\mathbf{v}\cdot \mathbf{n}) \, \ds\nolinebreak
		    - \int_{\gamma^-} (\gamma_0^-p) (-\mathbf{v} \cdot \mathbf{n}) \, \ds \\
		 &= \int_{\Omega} p \Div \mathbf{v} \, \dx \nolinebreak 
		     - \int_{\partial\Omega}p_D \mathbf{v} \cdot \mathbf{\nu} \, \ds \nolinebreak
		     - \int_{\gamma} t_f^{-1} (\mathbf{u}\cdot\mathbf{n})(\mathbf{v}\cdot\mathbf{n})\, \ds.
	\end{align*}
where the last step is derived by using the third equation of \eqref{BVP} and the definition of $\pjump$. Then, multiplying the first equation of \eqref{BVP} by $q\in \SpacepBVP$, we have:
    \begin{align*}
        \int_{\Omega} \Div \mathbf{u}q \, \dx = \int_{\Omega} f q \, dx.
    \end{align*}

We define $a: \Spaceu \times \Spaceu \xrightarrow{} \mathbb{R}$ and $b: \SpacepBVP \times \Spaceu \xrightarrow{} \mathbb{R}$ as bilinear functionals with the following forms:
	\begin{align*}
		a(\mathbf{u}, \mathbf{v}) &= (\perm^{-1}\mathbf{u},\mathbf{v} )_\Omega \nolinebreak
		 + ( t_f^{-1} \mathbf{u}\cdot\mathbf{n}, \mathbf{v}\cdot\mathbf{n} )_\gamma, \\
		b(p,\mathbf{v}) &= (p, \Div\mathbf{v} )_\Omega,
	\end{align*}
where $(\cdot,\cdot )_\Omega$ and $(\cdot,\cdot)_\gamma$ stand for the $L^2$ inner product on $\Omega$ and $\gamma$, respectively. Moreover, we define functionals $F:\SpacepBVP \xrightarrow{} \mathbb{R}$ and $G:\Spaceu \xrightarrow{} \mathbb{R}$ given by
	\begin{align*}
		F(q) = (f,q )_\Omega,\qquad G(\mathbf{v}) = \langle\mathbf{v}\cdot\mathbf{\nu} , p_D \rangle_{H^{-1/2}(\partial\Omega),H^{1/2}(\partial\Omega)}.
	\end{align*}
The mixed variational problem of \eqref{BVP} is then given by: Find $(\mathbf{u}, p)\in\mathcal{H}$ such that
   \begin{align}\label{vari_mix}
	\begin{split}
	a(\mathbf{u},\mathbf{v})-b(p,\mathbf{v}) &= -G(\mathbf{v}) \qquad \forall \mathbf{v} \in \Spaceu,\\
	b(q,\mathbf{u})  &= F(\mathbf{q}) \qquad \forall q\in\SpacepBVP.
	\end{split}
   \end{align}

\subsection{Derivation of the new model}
In this part, we first define $H^1(\Omega)$-approximations of pressure by considering the effects of the fault as source terms in expression of $-\Delta p$ and by assuming uniform normal direction of the fault. Then we derive the formulations for the approximate pressures by expressing the pressure jump along fault by the pressure derivatives so that the continuous Garlerkin finite element method could be used to solve the problem.

\subsubsection{Definition of approximations of pressure}

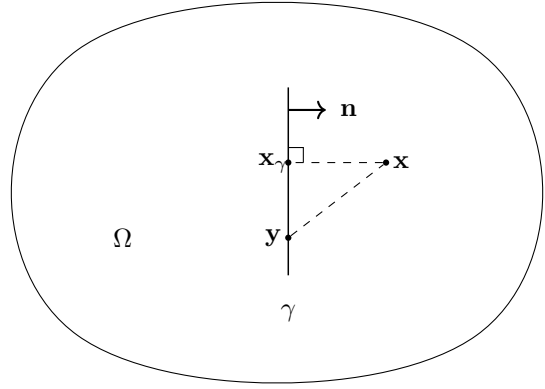
\begin{figure}[!t]
\begin{center}
  \begin{tikzpicture}
    \path[font={\tiny}]
        (11 , -.8)   coordinate (A)
        (5.7   , -.8)   coordinate (B)
        (5.7   , 3)   coordinate (C)
        (11   , 3)   coordinate (D)
    ;
    \draw plot [smooth cycle, tension=0.8] coordinates {(A) (B) (C) (D)};
    \draw[black, line width = 0.2mm] (8.5,0) to (8.5,2.5);
        \node[] at (6.3,0.5)   (a) {$\Omega$};
        
        \node[] at (8.3, 0.5) (A) {$\mathbf{y}$};
        \node[circle, fill, scale = 0.25] at (8.5,0.5)(a) {};
        
        \node[] at (10, 1.5) (B) {$\mathbf{x}$};
        \node[circle, fill, scale = 0.25] at (9.8,1.5)(b) {};
        \node[] at (8.5,-0.5) (e) {$\gamma$};
        \node[] at (8.3,1.5) (C) {$\mathbf{x_\gamma}$};
         \node[circle, fill, scale = 0.25] at (8.5,1.5)(c) {};
         
        \node[] at (9.3,2.2) (d) {$\mathbf{n}$};
    \draw[thick, ->] (8.5, 2.2) to (9,2.2);

    \draw[dashed] (a) -- (b)--(c)--cycle;
    \draw (8.7,1.5) -- (8.7,1.7) -- (8.5,1.7);
  \end{tikzpicture}
  \end{center}
  \caption{The domain $\Omega$ with an immersed fault $\gamma$: $\mathbf{y}$ and $\mathbf{x}$ are given points on $\gamma$ and in $\Omega$, $\mathbf{x}_\gamma$ is the projection of $\mathbf{x}$ on $\gamma$ along the normal direction $\mathbf{n}$.}
  \label{domain_point}
\end{figure}

By applying Green's Identities we can express $-\Delta p$ in the boundary value problem \eqref{BVP} by
    \begin{align}\label{L_p_D}
        -\Delta \pD\vari = f\vari - \int_\gamma\pjump\nabla \delta_{\mathbf{x}} \cdot \mathbf{n} \,\ds, \qquad \forall \mathbf{x} \in \Omega,
    \end{align}
where the derivation is presented in appendix \ref{appendices}. In addition, we assume the fault has uniform normal direction, i.e., fault as a point when $d=1$, a straight line when $d=2$, or a plane surface when $d=3$. Given a point $\mathbf{x} \in \Omega\backslash \gamma$, we denote the projection of $\mathbf{x}$ onto $\gamma$ by $\mathbf{x}_\gamma$. Using a coordinate basis $\{\mathbf{n}, \mathbf{\tau_1},\mathbf{\tau_2},\ldots,\mathbf{\tau_{d-1}} \}$, we can express $\mathbf{x} = (x_n, x_{\tau_1}, \ldots, x_{\tau_{d-1}} ) = (x_n, \mathbf{x_{\tau}} )$, and thus $\mathbf{x}_\gamma = (0, \mathbf{x_{\tau}} )$. A configuration of the domain is shown in \autoref{domain_point} where the fault can be expressed by $\gamma(\mathbf{y_\mathbf{\tau}}) = \left\{(0, \mathbf{y}_\mathbf{\tau} ): \mathbf{y}\in \gamma\right\} $. Thus, we can express the partial derivative of delta mass at $\mathbf{x}$ by
\begin{align*}
    \nabla \delta_{\mathbf{x}}(\mathbf{y}) \cdot \mathbf{n} \nolinebreak
     &= D_1 \delta_0(\mathbf{y}- \mathbf{x}) = D_1 \delta_0((-x_n, \mathbf{y_\mathbf{\tau}} - \mathbf{x_\mathbf{\tau}}))\nolinebreak
     = \delta^\prime_0 (-x_n)\delta_\mathbf{0} (\mathbf{y_\mathbf{\tau}} - \mathbf{x_\mathbf{\tau}}),
\end{align*}
where $\delta^\prime_0 (-x_n)$ represents the dipole distribution defined on one dimension, and $\delta_\mathbf{0} (\mathbf{y_\mathbf{\tau}} - \mathbf{x_\mathbf{\tau}})$ is the delta mass defined on $(d-1)$ dimension. By using this formulation, we can simplify the integration on $\gamma$ in \eqref{L_p_D} by
\begin{align*}
    \int_\gamma\pjump\nabla \delta_{\mathbf{x}} \cdot \mathbf{n} \, \ds &= \nolinebreak
    \int_{\gamma(\mathbf{y_\mathbf{\tau}})} \pjump(\mathbf{y_\mathbf{\tau}}) \delta^\prime_0 (-x_n)\delta_0 (\mathbf{y_\mathbf{\tau}} - \mathbf{x_\mathbf{\tau}}) \, \ds \\
     &=  -\delta^\prime_0 (x_n)\pjump (\mathbf{x_\mathbf{\tau}})\\
     &= -\delta^\prime_0 (x_n) t_f^{-1} u_{n,\gamma}(\mathbf{x_\mathbf{\tau}}),
\end{align*}
where $u_{n,\gamma}$ represents $(\mathbf{u}\cdot\mathbf{n})|_\gamma$ for simplicity and the last step is obtained by using the relationship between $\pjump$ and $u_{n,\gamma}$ implied by \eqref{BVP_fault}. By substituting the above formula into \eqref{L_p_D}, we have 
\begin{align}\label{p_tilde}
     -\Delta \pD\vari = f\vari + \delta^\prime_0 (x_n) t_f^{-1} u_{n,\gamma}(\mathbf{x_\mathbf{\tau}}) \qquad \forall \mathbf{x}=(x_n, \mathbf{x_{\mathbf{\tau}}}) \in \Omega. 
\end{align}

\theoremstyle{definition}
\begin{definition}
Suppose that $\left \{ \delta_\epsilon\right \}\subset \mathcal{D}(\Omega)$ is a sequence such that $\delta_\epsilon \xrightarrow{\mathcal{D}^\prime(\Omega)} \delta_0$ as $\epsilon \rightarrow 0^+$, we define the approximate sequence of pressures $\left \{ p_\epsilon \right \}$ as the solutions of the following boundary value problems:
\begin{align}\label{approx_BVP}
    \begin{cases}
    -\Delta p_\epsilon(\mathbf{x}) = f\vari + \dapprox (x_n) t_f^{-1} u_{n,\gamma}(\mathbf{x_\mathbf{\tau}}) \qquad &\forall \mathbf{x}=(x_n, \mathbf{x_{\mathbf{\tau}}}) \in \Omega,\\
    p = p_D \qquad &\text{on $\partial \Omega$.} 
    \end{cases}
\end{align}
\end{definition}

We conclude the properties of $p_\epsilon$ as the following theorem.
\begin{thm}
There is a unique $p_\epsilon\in H^1(\Omega)$ solving the problem of \eqref{approx_BVP} for each $\epsilon$. Moreover, the sequence of $p_\epsilon$ converges to $p$ in $L^2(\Omega)$ as $\epsilon \xrightarrow{} 0^+$.

\begin{myproof}
Define 
\begin{align*}
    \tilde{f}(\mathbf{x}) = f\vari + \dapprox (x_n) t_f^{-1} u_{n,\gamma}(\mathbf{x_\mathbf{\tau}}) \qquad &\forall \mathbf{x}=(x_n, \mathbf{x_{\mathbf{\tau}}}) \in \Omega.
\end{align*}
Since $\dapprox\in \mathcal{D}(\Omega)$ and $u_{n,\gamma}\in L^2(\Omega)$, then $\tilde{f}\in L^2(\Omega)$. Thus we expect $p_\epsilon \in H^1(\Omega)$ and the boundary vlaue problem \eqref{approx_BVP} is equivalent to find $p_\epsilon \in H^1_0(\Omega) + p_D$ such that: 
\begin{align}\label{vari_BVP_approx}
    a_\epsilon(p_\epsilon,q) = F_\epsilon(q) \qquad \forall q\in H^1_0(\Omega),
\end{align}
where $a_\epsilon:H^1(\Omega) \times H^1(\Omega) \xrightarrow{} \mathbb{R}$ is given by
\begin{align*}
    a_\epsilon(p_\epsilon, q) = (\nabla p_\epsilon, \nabla q)_{L^2(\Omega)},
\end{align*}
and $F_\epsilon:H^1_0(\Omega)\xrightarrow{} \mathbb{R}$ is given by
\begin{align*}
    F_\epsilon(q) = (\tilde{f}, q)_{L^2(\Omega)}.
\end{align*}
It is easy to notice that $a_\epsilon$ is continuous and coercive on $H^1(\Omega)$, therefore there is unique solution $p_\epsilon \in H_0^1(\Omega)$ by Lax-Milgram Theorem for each $\epsilon$.
Since $\delta_\epsilon^\prime \xrightarrow{\mathcal{D}^\prime(\Omega)} \delta_0^\prime$, $\Delta  p_\epsilon \xrightarrow{\mathcal{D}^\prime(\Omega)} \Delta \pD$. Take any $\phi\in \mathcal{D}(\Omega)$, we have
\begin{align*}
    \langle p_\epsilon, \Delta\phi \rangle_\Omega = \langle \Delta p_\epsilon, \phi\rangle_\Omega + \nolinebreak
    \langle p_D, \nabla \phi \cdot \mathbf{\nu} \rangle_{\partial \Omega} - \langle \nabla p_\epsilon \cdot \mathbf{\nu}, \phi \rangle_{\partial \Omega}\\
    \longrightarrow \langle \Delta \pD, \phi\rangle_\Omega + \nolinebreak
    \langle p_D, \nabla \phi \cdot \mathbf{\nu} \rangle_{\partial \Omega} - \langle \nabla \pD \cdot \mathbf{\nu}, \phi \rangle_{\partial \Omega} = \langle \pD, \Delta \phi \rangle_\Omega.
\end{align*}
Thus we can conclude that $p_\epsilon \xrightarrow{\mathcal{D}^\prime(\Omega)} \pD$, and thus $p_\epsilon \xrightarrow{L^2(\Omega)} \pD$ by Lebesgue lemma. 
\end{myproof}
\end{thm}

Note that the normal component of velocity along the fault ($u_{n, \gamma}$ in \eqref{approx_BVP}) is unknown in our problem, so we can not directly solve \eqref{approx_BVP} to get $p_\epsilon$ by continuous Galerkin method. But if we can express $u_{n,\gamma}$ by derivatives of $p$, we can formulate differential equations about $p_\epsilon$.

\subsubsection{Formulation for \texorpdfstring{$d=1$}{Lg}}
\paragraph{Model derivation.} Consider a one-dimensional domain $\Omega = (0,L)$ and a fault as a point at $x_\gamma \in \Omega$. We additionally assume that there is no external source such that $f = 0$, and that the boundary conditions imposed are $p(0) = p_0, \, p(L) = p_L$. Then the original boundary value problem \eqref{BVP} has pressure solution as
\begin{align}\label{eq_d1_anp}
    p_c(x) = \begin{cases} 
    p_0 - \frac{p_0-p_L}{t_f^{-1} + L}x  &0 \leq x < x_\gamma,\\
    p_L + \frac{p_0-p_L}{t_f^{-1} + L}(L-x) &x_\gamma< x \leq L,
    \end{cases}
\end{align}
and velocity solution as
\begin{align}\label{eq_d1_anv}
    u_c(x) = u_{n,\gamma} = \frac{p_0 - p_L}{t_f^{-1} + L} \qquad 0\leqslant x \leqslant L.
\end{align}
\eqref{approx_BVP} for the approximate pressure $p_\epsilon$ in one dimension could be expressed by
\begin{align}\label{eq_d1_approx}
    -\diff[2]{p_\epsilon}{x} =  \deltaeps^\prime \frac{u_{n,\gamma}}{t_f} \qquad & 0 < x < L.
\end{align}
In \autoref{1d_illustration}, we plotted $p_c$ and several $p_\epsilon$ under different values of $\epsilon$ as illustration. If $u_{n,\gamma}$ is given, $p_\epsilon$ could be expressed explicitly in analytical solution $p_c$ with the pressure jump along the fault, i.e., $\pjumpinu$ by
\begin{align}\label{eq_ana_p}
    p_\epsilon^c(x) = 
    \begin{cases}
        p_c(x) - H_\epsilon (x) \pjumpinu \qquad &0\leqslant x \leqslant x_\gamma,\\
        p_c(x) - (H_\epsilon(x) -1) \pjumpinu \qquad &x_\gamma \leqslant x \leqslant L,
    \end{cases}
\end{align}
where 
\begin{align*}
    H_\epsilon(x) =\int_{0}^{x} \deltaeps(t) \, dt.
\end{align*}
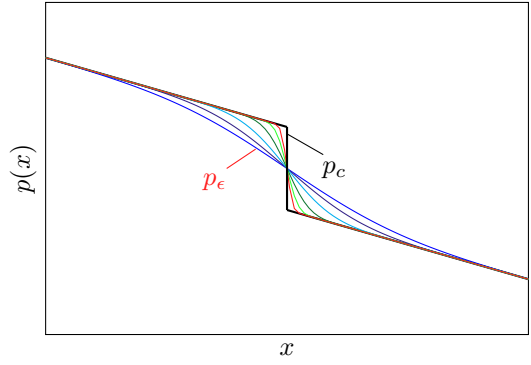
\begin{figure}[!t]
\centering
\begin{tikzpicture}[align=right,yscale=1.,xscale=1.]
\begin{groupplot}[group style={group size=2 by 1,horizontal sep=2cm},
height=6cm, width=8cm,
ticks=none,
xmin=-1, xmax=1,
ymin = -1.2, ymax = 1.2]
\nextgroupplot[xlabel=$x$, ylabel =$p(x)$, xlabel near ticks, ylabel near ticks]
    \addplot[black, thick, domain=-1:0.] {-0.5*x +.3};
    \addplot[black, thick, domain=0:1.,samples=100] {-0.5*x -.3};
    \addplot[blue, domain=-1.:1.,samples=100]  {0.3*erf(-2.*x)-0.5*x};
    \addplot[c4, domain=-1.:1.,samples=100]  {0.3*erf(-3.*x)-0.5*x};
    \addplot[cyan, domain=-1.:1.,samples=100] {0.3*erf(-5.*x)-0.5*x};
    \addplot[c3, domain=-1.:1.,samples=100] {0.3*erf(-10.*x)-0.5*x};
    \addplot[green!80, domain=-1.:1.,samples=100]  {0.3*erf(-20.*x)-0.5*x};
    \addplot[red!90, domain=-1.:1.,samples=100]  {0.3*erf(-40.*x)-0.5*x};
    \draw[red!90] (axis cs:-0.13,0.14) -- node[below left = 0.1cm] {$p_\epsilon$} (axis cs:-0.25,0);
    \draw[black] (axis cs:0.0,0.25) -- node[below right= 0.15cm] {$p_c$} (axis cs:0.15,0.1);
    \draw[black,thick] (axis cs:0.,-0.3) -- (axis cs:0.,0.3);
\end{groupplot}
\end{tikzpicture}
\caption{An illustration of the relationship between correct pressure solution $p_c$ and the approximated solutions $p_\epsilon$ in one-dimensional space.}
\label{1d_illustration}
\end{figure}

Since $f=0$, it is natural to expect that the solution of the above equation is twice differentiable, i.e., $p_\epsilon \in C^2([0,L])$. By integrating \eqref{approx_BVP}, one can get
\begin{align}\label{eq_d1_2}
    -\diff{p_\epsilon}{x} = M + \deltaeps \frac{u_{n,\gamma}}{t_f} ,
\end{align}
where $M\in\mathbb{R}$ is a constant. Integrate the above equation from $x=0$ to $x=L$, we have
\begin{align*}
    \int_{0}^{L} -\diff{p_\epsilon}{x} \, dx = p_0 - p_L = ML + \frac{u_{n,\gamma}}{t_f}.
\end{align*}
But the analytical solution of pressure shows that
\begin{align*}
    p_0 - p_L = u_{n,\gamma}L + \pjump= u_{n,\gamma}L + \frac{u_{n,\gamma}}{t_f}.
\end{align*}
Thus we can conclude that $M = u_{n,\gamma}$. Substituting $M$ into \eqref{eq_d1_2}, we can express $u_{n,\gamma}$ as
\begin{align}\label{eq_d1_u}
    u_{n,\gamma} = -(1+ \frac{\deltaeps}{t_f})^{-1} \diff{p_\epsilon}{x}.
\end{align}
Using this expression of $u_{n,\gamma}$, we convert \eqref{eq_d1_approx} into the following new equation:
\begin{align}\label{eq_d1_p}
    -\diff[2]{p_\epsilon}{x} +  \frac{\deltaeps^\prime}{t_f+\deltaeps} \diff{p_\epsilon}{x}=0 \qquad 0< x < L.
\end{align}

\begin{figure}[!t]
	\centering
	\begin{tikzpicture}[every plot/.append style={thick}]
	\begin{groupplot}[
        footnotesize, 
        width=6cm,
        height=5cm,
        group style={
                     group name=1d, group size=3 by 2,
                     horizontal sep=1.0cm, vertical sep=1.0cm, 
                     },
        xlabel near ticks, ylabel near ticks,
        yticklabel style={/pgf/number format/fixed}
      ]
     \nextgroupplot[ legend cell align={left},
   						legend style={at={(axis cs:7,0.55)},anchor=south west, draw = none,line width=1.0pt},
                        xmin = 0, xmax = 10,
                        xtick ={0,2,...,10},
                        xlabel=$x$,
                        ylabel=$p$,
                        ytick ={0,0.2,...,1},
                        ymin = 0, ymax = 1]

                        \addplot[red] table [x=x ,y=p_exact] {plot/1D_1.0.txt};
             		    \addplot[blue!70,dotted] table [x=x ,y=p] {plot/1D_1.0.txt};
             		    \addplot[c3,dashed] table [x=x ,y=p] {plot/1D_1.0.txt};

            \addlegendimage{/pgfplots/refstyle=plot_one}
            \addlegendentry{\scriptsize{$p_c$}}
            \addlegendentry{\scriptsize{$p_\epsilon^c$}}
            \addlegendentry{\scriptsize{$p_\epsilon$}}
            	\label{plot_one} 

         \nextgroupplot[ legend cell align={left},
   						legend style={at={(axis cs:7,0.55)},anchor=south west, draw = none,line width=1.0pt},
                        xmin = 0, xmax = 10,
                        xtick ={0,2,...,10},
                        xlabel=$x$,
                        ytick ={0,0.2,...,1},
                        ymin = 0, ymax = 1]

                        \addplot[red] table [x=x ,y=p_exact] {plot/1D_0.5.txt};
             		    \addplot[blue!70,dotted] table [x=x ,y=p] {plot/1D_0.5.txt};
             		    \addplot[c3,dashed] table [x=x ,y=p] {plot/1D_0.5.txt};

            \addlegendimage{/pgfplots/refstyle=plot_one}
            \addlegendentry{\scriptsize{$p_c$}}
            \addlegendentry{\scriptsize{$p_\epsilon^c$}}
            \addlegendentry{\scriptsize{$p_\epsilon$}}
            	\label{plot_one}

             \nextgroupplot[ legend cell align={left},
   						legend style={at={(axis cs:7,0.55)},anchor=south west, draw = none,line width=1.0pt},
                        xmin = 0, xmax = 10,
                        xtick ={0,2,...,10},
                        xlabel=$x$,
                        ytick ={0,0.2,...,1},
                        ymin = 0, ymax = 1]

                        \addplot[red] table [x=x ,y=p_exact] {plot/1D_0.01.txt};
             		    \addplot[blue!70,dotted] table [x=x ,y=p] {plot/1D_0.01.txt};
             		    \addplot[c3, dashed] table [x=x ,y=p] {plot/1D_0.01.txt};

            \addlegendimage{/pgfplots/refstyle=plot_one}
            \addlegendentry{\scriptsize{$p_c$}}
            \addlegendentry{\scriptsize{$p_\epsilon^c$}}
            \addlegendentry{\scriptsize{$p_\epsilon$}}
            	\label{plot_one} 
        
    \nextgroupplot[ legend cell align={left},
   						legend style={at={(axis cs:7,0.035)},anchor=south west, draw = none,line width=1.0pt},
                        xmin = 0, xmax = 10,
                        xtick ={0,2,...,10},
                        xlabel=$x$,
                        ylabel=$u$,
                        ytick ={0,0.02,...,0.1},
                        ymin = 0, ymax = 0.1]

                        \addplot[red] table [x=x ,y=u_exact] {plot/1D_1.0.txt};
             		    \addplot[c3,dashed] table [x=x ,y=u] {plot/1D_1.0.txt};

            \addlegendimage{/pgfplots/refstyle=plot_one}
            \addlegendentry{\scriptsize{$u_c$}}
            \addlegendentry{\scriptsize{$u_\epsilon$}}
            	\label{plot_one} 
            	
        \nextgroupplot[ legend cell align={left},
   						legend style={at={(axis cs:7,0.035)},anchor=south west, draw = none,line width=1.0pt},
                        xmin = 0, xmax = 10,
                        xtick ={0,2,...,10},
                        xlabel=$x$,
                        ytick ={0,0.02,...,0.1},
                        ymin = 0, ymax = 0.1]

                        \addplot[red] table [x=x ,y=u_exact] {plot/1D_0.5.txt};
             		    \addplot[c3,dashed] table [x=x ,y=u] {plot/1D_0.5.txt};

            \addlegendimage{/pgfplots/refstyle=plot_one}
            \addlegendentry{\scriptsize{$u_c$}}
            \addlegendentry{\scriptsize{$u_\epsilon$}}
            	\label{plot_one}
            	
        \nextgroupplot[ legend cell align={left},
   						legend style={at={(axis cs:7,0.035)},anchor=south west, draw = none,line width=1.0pt},
                        xmin = 0, xmax = 10,
                        xtick ={0,2,...,10},
                        xlabel=$x$,
                        ytick ={0,0.02,...,0.10},
                        ymin = 0, ymax = 0.10]

                        \addplot[red] table [x=x ,y=u_exact] {plot/1D_0.01.txt};
             		    \addplot[c3,dashed] table [x=x ,y=u] {plot/1D_0.01.txt};

            \addlegendimage{/pgfplots/refstyle=plot_one}
            \addlegendentry{\scriptsize{$u_c$}}
            \addlegendentry{\scriptsize{$u_\epsilon$}}
            	\label{plot_one}

	\end{groupplot}
	\node[above = .5cm of 1d c1r1.north] {$\epsilon = 1.0$};
    \node[above = .5cm of 1d c2r1.north] {$\epsilon = 0.5$};
    \node[above = .5cm of 1d c3r1.north] {$\epsilon = 0.01$};

	\end{tikzpicture}
	\caption{Pressure and velocity solutions in the 1D test.}
	\label{1d_plot}
\end{figure}
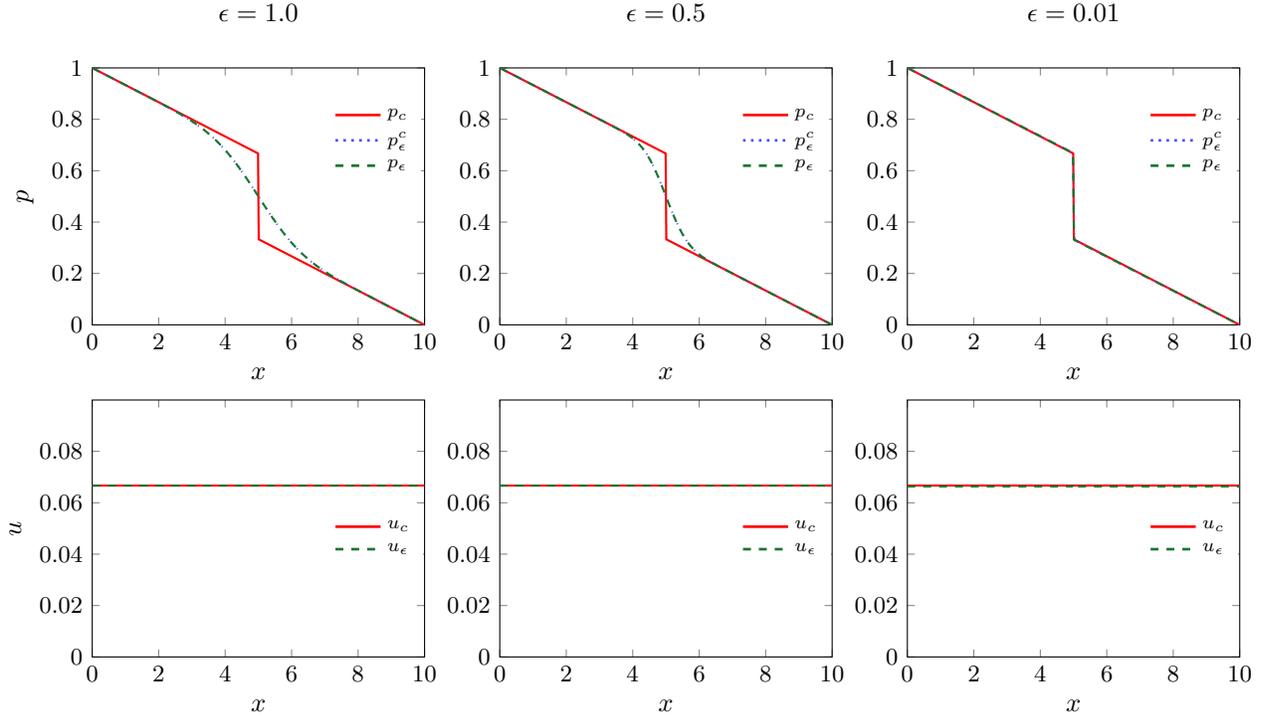

\paragraph{Numerical test.} To show the equivalence between \eqref{eq_d1_approx} and \eqref{eq_d1_p}, we consider $\Omega = (0,10)$, a fault at $x_\gamma=5$ with transmissibility $t_f = 0.2$. Choose $p_0 = 1, \,p_L=0$ as boundary condition, and use Gaussian functions to approximate delta distribution:
\begin{align*}
    \deltaeps(x)  = \frac{1}{\sqrt{2\pi}\epsilon} \exp{\Big(-\frac{(x-x_\gamma)^2}{2\epsilon^2}\Big)}.
\end{align*}
Though Gaussian functions do not have compact support in $\Omega$, we can assume that $\deltaeps \rightarrow \delta_0(x_n)$ and$\int_{\Omega} \deltaeps(x) dx = 1$ since $\epsilon$ we considered here are much smaller than the range of the domain. \eqref{eq_d1_p} is solved by continuous Galerkin method using Lagrange basis function with degree of 1.

We denote the analytical solutions of \eqref{eq_d1_anp} and \eqref{eq_d1_anv} by $p_c$ and $u_c$, and approximate pressure solution of \eqref{eq_d1_approx} by $p_{\epsilon}^c$, and the numerical solutions of our derived model \eqref{eq_d1_p} and \eqref{eq_d1_u} by $p_\epsilon$ and $u_\epsilon$. In the experiments, we consider three values of $\epsilon = 1.0, \, 0.5, \, \text{and } 0.01$. Comparison results of pressure and velocity solutions between different models are presented in \autoref{1d_plot}. From the pressure solutions, $p_\epsilon$ match $p_\epsilon^c$ under different values of $\epsilon$. Moreover, As $\epsilon$ becomes smaller, the approximate continuous pressure is getting closer to the discontinuous analytical solution. Furthermore, the derived model for velocity can give correct velocity solution under different values of $\epsilon$.

\subsubsection{Formulation for \texorpdfstring{$d\geqslant 2$}{Lg}}
\paragraph{Strong formulation.} \eqref{approx_BVP} for the approximated pressure $p_\epsilon$ which is denoted as $p$ here for simplicity can be expressed by
\begin{align}\label{eq_d2_approxp}
    -\diffp[2]{p_\epsilon}{\mathbf{n}}- \sum_{i=1}^{d-1} \diffp[2]{p_\epsilon}{\mathbf{\tau_i}} = \deltaeps^\prime \pjumpinu.
\end{align}
in reference to the coordinates of the fault $\{\mathbf{n}, \mathbf{\tau_1},\mathbf{\tau_2},\ldots,\mathbf{\tau_{d-1}} \}$.

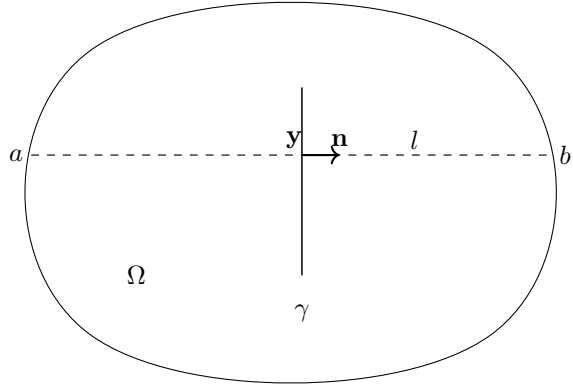
\begin{figure}[!t]
\begin{center}
  \begin{tikzpicture}
    \path[font={\tiny}]
        (11 , -.8)   coordinate (A)
        (5.7   , -.8)   coordinate (B)
        (5.7   , 3)   coordinate (C)
        (11   , 3)   coordinate (D)
    ;
    \draw plot [smooth cycle, tension=0.8] coordinates {(A) (B) (C) (D)};
    \draw[black, line width = 0.2mm] (8.5,0) to (8.5,2.5);
        \node[] at (6.3,0)   (a) {$\Omega$};
        \node[] at (8.5,-0.5)   (a) {$\gamma$};
        \draw[dashed, black, line width = 0.15mm] (4.9,1.6) to (11.8,1.6);
        \node[] at (10,1.8)   (b) {$l$};
        \node[] at (4.7,1.6)   (c) {$a$};
        \node[] at (12,1.6)   (d) {$b$};
        \node[] at (8.4,1.8)   (e) {$\mathbf{y}$};
        \node[] at (9,1.8) (c) {$\mathbf{n}$};
        \draw[thick, ->] (8.5, 1.6) to (9,1.6);
  \end{tikzpicture}
  \end{center}
  \caption{The domain $\Omega$ with an immersed fault $\gamma$: for each point $\mathbf{y}$ of $\gamma$, we consider a line $l\in\Omega$ perpendicularly intersects with $\gamma $ at $\mathbf{y}$ ($a$ and $b$ are points of $l$ intersecting $\partial \Omega$.}
  \label{domain_d2}
\end{figure}

We denote the pressure solution of the original boundary value problem \eqref{BVP} by $ p_c$. Compared with the linear pressure profile when $d=1$, the pressure in general is nonlinear when $d\geqslant2$. But we can decompose the pressure into two parts such that
\begin{align*}
    p_c = p_{c,l} + g_c,
\end{align*}
where $p_{c,l}$ represents the "linear part "of $p_c$ (in $\mathbf{n}$ with respect to the fault) such that
\begin{align*}
    p_{c,l}|_{\gamma^+}=p_c|_{\gamma^+}, \qquad p_{c,l}|_{\gamma^-}=p_c|_{\gamma^-}, \qquad 
\end{align*}
and
\begin{align*}
    -\nabla p_{c,l}\cdot \mathbf{n}= -\nabla p_{c}\cdot \mathbf{n}|_\gamma = u_{n,\gamma}.
\end{align*}
Thus, the "nonlinear part" $g_c$ must satisfy the following conditions: 
\begin{align*}
    g_c|_{\gamma^+} =g_c|_{\gamma^-} =0, \qquad     (\nabla g_c\cdot \mathbf{n})|_{\gamma}= 0.
\end{align*}
From the previous construction, it is clearly that such decomposition is unique for a given $p_c$.\par

\begin{figure}[!t]
\centering
\begin{tikzpicture}[align=right,yscale=1.,xscale=1.,every plot/.append style={thick}]
\begin{groupplot}[group style={group name = 2d, group size=2 by 1,horizontal sep=2cm},
height=6cm, width=8cm,
ticks=none,
xmin=-1, xmax=1,
ymin = -1.2, ymax = 1.2]
\nextgroupplot[xlabel=$x_n$, ylabel =$p(x_n)$, xlabel near ticks, ylabel near ticks]
    \addplot[blue!70, domain=-1:0,mark=none]{0.3*(x^2.) - 0.5*x + .3};
    \addplot[gray, domain=-1:0.] {-0.5*x +.3};
      \addplot[gray, domain=0:1.,samples=30] {-0.5*x -.3};
    \addplot[blue!70, domain=0:1.,samples=30] {-0.3*(x^2.) - 0.5*x - .3};
    \draw[thick,gray] (axis cs:0.,-0.3) -- (axis cs:0.,0.3);
     \draw[black,<->,thick] (axis cs:-0.8,0.892) -- node[left] {$g_c$} (axis cs:-0.8,0.7);
     \draw[black] (axis cs:-0.6,0.59) -- node[below right = 0.1cm] {$p_{c,l}$} (axis cs:-0.5,0.4);
    \draw[blue!70] (axis cs:-0.4,0.548) -- node[above right = 0.1cm] {$p_c$} (axis cs:-0.25,0.7);
    \draw (axis cs:0.,-0.3) -- (axis cs:0.1,-0.3);
    \draw (axis cs:0.,0.3) -- (axis cs:0.1,0.3);
    \draw[dashed,<->,thick] (axis cs:0.05,-0.3) -- node[right=0.1cm]{$\pjump$} (axis cs:0.05, 0.3);

\nextgroupplot[xlabel=$x_n$, ylabel =$p(x_n)$, xlabel near ticks, ylabel near ticks]
    \addplot[blue!70, dashed, domain=-1:0,mark=none]{0.3*(x^2.) - 0.5*x + .3};
    \addplot[gray, dashed, domain=-1:0.] {-0.5*x +.3};
      \addplot[gray, dashed, domain=0:1.,samples=30] {-0.5*x -.3};
     \addplot[blue!70,dashed, domain=0:1.,samples=30] {-0.3*(x^2.) - 0.5*x - .3};
    \addplot[red!90, domain=-1.:0.,samples=30] {0.3*erf(-3.*x)-0.5*x};
     \addplot[red!90, domain=0.:1.,samples=30] {0.3*erf(-3.*x)-0.5*x};
      \addplot[green!70!black, domain=-1.:0.,samples=30] {0.3*erf(-3.*x)-0.5*x +0.3*(x^2)};
    \addplot[green!70!black, domain=0.:1.,samples=30] {0.3*erf(-3.*x)-0.5*x -0.3*(x^2)};
    \draw[red!90] (axis cs:-0.6,0.59) -- node[below right = 0.1cm] {$p_l$} (axis cs:-0.5,0.4);
         \draw[black,thick,<->] (axis cs:-0.8,0.892) -- node[left] {$g$} (axis cs:-0.8,0.7);
    \draw[green!70!black] (axis cs:-0.4,0.548) -- node[above right = 0.1cm] {$p$} (axis cs:-0.25,0.7);
     \draw[dashed,gray,thick] (axis cs:0.,-0.3) -- (axis cs:0.,0.3);
\end{groupplot}
    \node[below = .5cm of 2d c1r1.south] {(a)};
    \node[below = .5cm of 2d c2r1.south] {(b)};
\end{tikzpicture}
\caption{An illustration of the functions defined in the case of multi-dimensions, }
\label{2d_online}
\end{figure}
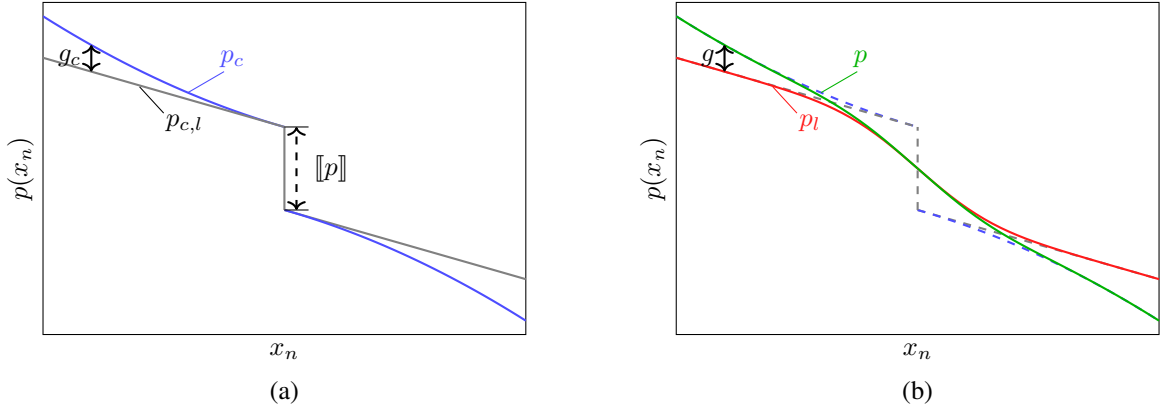

Take a point $\mathbf{y}=(0,\mathbf{y_\tau})\in\gamma$, we consider the part of line inside $\Omega$ that is perpendicularly intersecting $\gamma$ at $\mathbf{y}$, i.e., $l = \{ (x_n,\mathbf{y_\tau})\in\Omega:x_a\leqslant x_n \leqslant x_b\}$ in \autoref{domain_d2}. Thus $p_{c,l}(\cdot, \mathbf{y_\tau})$ satisfies the boundary value problem \eqref{BVP} in one dimension. The relationship between $p_c(\cdot, \mathbf{y_\tau)}$ and $p_{c,l}(\cdot, \mathbf{y_\tau})$ is sketched in \autoref{2d_online}(a). From the previous discussion in one dimension, there exists $\eta_{\mathbf{y_\tau}} \in C^2([x_a,x_b])$ as an approximation to $p_{c,l}(\cdot, \mathbf{y_\tau})$ such that,
\begin{align*}
    -\diff[2]{\eta_{\mathbf{y_\tau}}}{x_n} = \deltaeps^\prime \frac{u_{n,\gamma}}{t_f}, \qquad \nolinebreak
    \eta_{\mathbf{y_\tau}}(x_a) = p_{c,l}(x_a, \mathbf{y_\tau}),\qquad
    \eta_{\mathbf{y_\tau}}(x_b) = p_{c,l}(x_b, \mathbf{y_\tau}).
\end{align*}
Since $\mathbf{y}$ is a random point from $\gamma$, we can define $p_{\epsilon, l}$ as an approximation to $p_{c,l}$ such that $p_{\epsilon,l}(\cdot, \mathbf{y_\tau}) = \eta_{\mathbf{y_\tau}}(\cdot) $. Then,
\begin{align*}
    -\pdv[2]{p_{\epsilon,l}}{\mathbf{n}} = \deltaeps^\prime \frac{u_{n,\gamma}}{t_f}.
\end{align*}

Similar to the decomposition of $p_c$, we can decompose $p_\epsilon$ into $p_{\epsilon,l}$ and $g_\epsilon$ such that
\begin{align*}
    p_\epsilon = p_{\epsilon,l} + g_\epsilon.
\end{align*}
We sketched $p_{\epsilon, l}$, $g_\epsilon$ and $p_\epsilon$ in \autoref{2d_online}(b). Moreover, since $p_\epsilon \rightarrow p_c$ and $p_{\epsilon, l} \rightarrow p_{c,l}$ as $\epsilon \rightarrow 0^+$, we have $g_\epsilon \rightarrow g_c$ as $\epsilon \rightarrow 0^+$. Then, the first-order and second-order derivatives of $p_\epsilon$ in the normal direction can be expressed by
\begin{align}\label{eq_d2_dp}
        -\pdv{p_\epsilon}{\mathbf{n}} = -\pdv{p_{\epsilon ,l}}{\mathbf{n}}-\pdv{g_\epsilon}{\mathbf{n}} \nolinebreak
        = (t_f + \deltaeps) \pjumpinu - \pdv{g_\epsilon}{\mathbf{n}},
\end{align}
and
\begin{align*}
        -\pdv[2]{p_\epsilon}{\mathbf{n}} = \deltaeps ^\prime \pjumpinu - \pdv[2]{g_\epsilon}{\mathbf{n}}.
\end{align*}
Substitute the above equation into \eqref{eq_d2_approxp}, we have
\begin{align*}
    -\pdv[2]{g_\epsilon}{\mathbf{n}} = \sum_{i=1}^{d-1} \pdv[2]{p_\epsilon}{\mathbf{\tau}_i}.
\end{align*}
We assume that $\partial g_\epsilon/\partial \mathbf{n} = 0$ on $\gamma$, and approximate $-\partial^2 g_\epsilon/\partial \mathbf{n}^2$ by applying first-order Taylor expansion:
\begin{align*}
    -\pdv{g_\epsilon}{\mathbf{n}} (\mathbf{x}) =-x_n \pdv[2]{g_\epsilon}{\mathbf{n}} (\mathbf{x})  = x_n \sum_{i=1}^{d-1}  \pdv[2]{p_\epsilon}{\mathbf{\tau}_i} (\mathbf{x}) \qquad \forall \mathbf{x} = (x_n, \mathbf{x_\mathbf{\tau}})\in\Omega.
\end{align*}
Substitute the above equation to \eqref{eq_d2_dp}, we have the expression for the pressure jump given by
\begin{align}
    \pjumpinu = - (t_f + \deltaeps)^{-1} (\pdv{p_\epsilon}{\mathbf{n}} + x_n \sum_{i=1}^{d-1} \pdv[2]{p_\epsilon}{\mathbf{\tau}_i}).
\end{align}
Finally, by substituting the above expression of $\pjumpinu$ into \eqref{eq_d2_approxp}, we can derive a new equation for $p_\epsilon$ given by
\begin{align}\label{eq_d2_strong}
    -\Delta p_\epsilon = f - \frac{\deltaeps ^\prime}{t_f + \deltaeps}(\pdv{p_\epsilon}{\mathbf{n}} + x_n \sum_{i=1}^{d-1} \pdv[2]{p_\epsilon}{\mathbf{\tau}_i}) \qquad \forall \mathbf{x} = (x_n, \mathbf{x_\mathbf{\tau}})\in\Omega.
\end{align}

\paragraph{Weak formulation.} In order to get the weak formulation of the derived formulation \eqref{eq_d2_strong}, we first define functions $\dipole, \Dipole \in \mathcal{D}(\Omega)$ by 
\begin{align*}
    \dipole = \frac{\deltaeps^\prime}{t_f + \deltaeps},
    \qquad \Dipole = \frac{\deltaeps^\prime \, x_n}{t_f + \deltaeps}, \qquad \forall \mathbf{x} = (x_n, \mathbf{x_\mathbf{\tau}})\in\Omega.
\end{align*}
Thus, \eqref{eq_d2_strong} can be expressed by
\begin{align*}
    -\Delta p_\epsilon = f - \dipole \pdv{p_\epsilon}{\mathbf{n}} -  \sum_{i=1}^{d-1}\Dipole \pdv[2]{p_\epsilon}{\mathbf{\tau}_i}
\end{align*}

Take any $w\in H_0^1(\Omega)$, we have the weak formulation as
\begin{align}\label{weak_1}
    \Big(\nabla p_\epsilon, \nabla w\Big)_\Omega + \Big(\dipole \pdv{p_\epsilon}{\mathbf{n}}, w\Big)_\Omega + \nolinebreak
     \sum_{i=1}^{d-1} \Big(\Dipole \pdv[2]{p_\epsilon}{\mathbf{\tau}_i},w \Big)_\Omega = \Big(f,w\Big)_\Omega.
\end{align}
For each $i=1,2,...,d-1$, we define a vector $\mathbf{p}_i = \begin{bmatrix}0 &... &\pdv{p_\epsilon}{\mathbf{\tau}_i} &... &0\end{bmatrix}^T$ where $\pdv{p}{\mathbf{\tau}_i}$ is the $i$th component of $\mathbf{p}_i$. Then we can integrate by parts for the term related to $\pdv[2]{p_\epsilon}{\mathbf{\tau}_i}$ in \eqref{weak_1}
\begin{align*}
    \Big(\Dipole \pdv[2]{p_\epsilon}{\mathbf{\tau}_i}, w \Big)_\Omega &= \int_\Omega (\Dipole w) \nabla \cdot \mathbf{p}_i  \, \dx\nolinebreak
    = -\int_\Omega \mathbf{p_\epsilon}_i \nabla \cdot (\Dipole w) \, \dx + \int_{\partial \Omega} (\Dipole w) \mathbf{p}_i \cdot \mathbf{\nu} \, \dx \\
     &= -\Big(\pdv{\Dipole}{\mathbf{\tau}_i} \pdv{p_\epsilon}{\mathbf{\tau}_i}, w\Big)_\Omega - \Big(\Dipole \pdv{p_\epsilon}{\mathbf{\tau}_i}, \pdv{w}{\mathbf{\tau}_i}\Big)_\Omega \qquad i = 1,2,...,d-1.
\end{align*}

Define bilinear functional $B:H^1(\Omega) \times H^1(\Omega) \rightarrow \mathbb{R}$ by
\begin{align*}
    B_\epsilon (p_\epsilon, w) = \Big(\nabla p_\epsilon, \nabla w\Big)_\Omega + \Big(\dipole \pdv{p_\epsilon}{\mathbf{n}}, w\Big)_\Omega \nolinebreak
    - \sum_{i=1}^{d-1}\Big(\pdv{\Dipole}{\mathbf{\tau}_i} \pdv{p_\epsilon}{\mathbf{\tau}_i}, w\Big)_\Omega \nolinebreak
    - \sum_{i=1}^{d-1} \Big(\Dipole \pdv{p_\epsilon}{\mathbf{\tau}_i}, \pdv{w}{\mathbf{\tau}_i}\Big)_\Omega.
\end{align*}
Then the variational problem is to find $p_\epsilon\in H_0^1(\Omega) + p_D$ such that:
\begin{align}\label{weak_2}
    B_\epsilon(p_\epsilon,w) = F(w) \qquad \forall w \in H_0^1(\Omega).
\end{align}

\section{Implementation} \label{implementation}

\begin{figure}[!t]
\begin{center}
  \begin{tikzpicture}
     \path[font={\tiny}]
        (11 , -.8)   coordinate (A)
        (5.7   , -.8)   coordinate (B)
        (5.7   , 3)   coordinate (C)
        (11   , 3)   coordinate (D)
    ;
    \draw plot [smooth cycle, tension=0.8] coordinates {(A) (B) (C) (D)};
    \draw[black, line width = 0.2mm] (8.5,0) to (8.5,2.5);
    \node[] at (11,0.5)   (a) {$\Omega$};
    \node[] at (4.5,0.5)   (a) {$\partial\Omega$};
    \node[] at (8., 0)   (b) {$\Omega_\epsilon^s$};
    \node[] at (8.25, 1.5)   (b) {$\gamma$};
    \node[] at (7.1, 1)   (b) {$\Gamma_{\epsilon}^{s,+}$};
    \node[] at (9.9, 1)   (b) {$\Gamma_{\epsilon}^{s,-}$};
    \draw[blue!70] (7.5, -0.2) -- (9.5, -0.2) -- (9.5, 2.7) -- (7.5,2.7) -- (7.5, -0.2) ;
    \fill[pattern=flexible north east lines,pattern line width=0.2pt,pattern color=blue] (7.5,-0.2) rectangle (9.5,2.7);
  \end{tikzpicture}
  \end{center}
  \caption{The domain $\Omega$ with an immersed fault $\gamma$.}
  \label{domain_sub}
\end{figure}
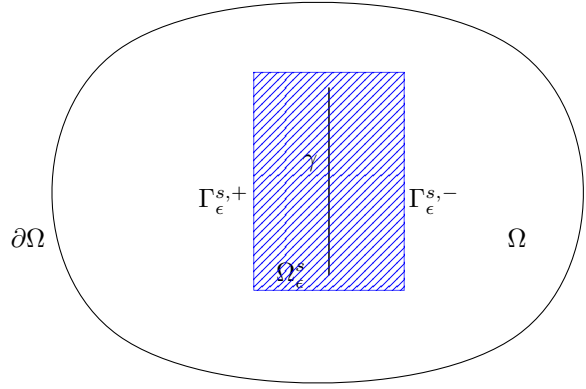

\subsection{New method procedures}
Based on the formulations we obtained for the approximations of pressure, we formulate a new method to obtain the pressure and velocity solutions. We first solve the variational problem of \autoref{weak_2} to get approximate pressure $p_\epsilon$ at given value of $\epsilon$. Then we solve $\mathbf{u}_\epsilon$ by Darcy's law:
\begin{align}\label{u_eps}
    \mathbf{u}_\epsilon = - \nabla p_\epsilon.
\end{align}

Note that $p_\epsilon$ approximates $p$ in the sense that it smooths the singularity of $p$ occurs at the fault. Moreover, the normal component of velocity $\mathbf{u}$ on the fault in the original problem \eqref{BVP} is not defined by \eqref{u_eps}, but by the trace of $\mathbf{u\cdot n}$. Thus, we can expect that $\mathbf{u}_\epsilon$ may have large difference around the fault compared to the correct velocity solution $\mathbf{u}$. Therefore, in order to get sufficiently accurate solutions, we additionally solve problem associated with a small domain around fault to modify $p_\epsilon$ and $\mathbf{u}_\epsilon$.

Suppose that our defined function $\text{supp}(\delta_\epsilon) \subset \Omega_\epsilon^s$, which is an open subset of $\Omega$ for given $\epsilon$, we can expect that $p_\epsilon$ should be the same as $p$ in $\Omega \backslash \Omega_\epsilon^s$, so as $\mathbf{u}_\epsilon$. This inspires that we only need to modify our solutions in a subdomain around fault, and our solutions are accurate outside the subdomain. Assume that $\Omega_\epsilon^s$ is d-orthotope, i.e. rectangle when $d=2$ and cuboid when $d=3$, such that it has part of boundary is parallel to $\gamma$, denote the part on the side of $\gamma^+$ by $\Gamma_{\epsilon}^{s,+}$, and the other part by $\Gamma_{\epsilon}^{s,-}$ (see \autoref{domain_sub}). Define the following boundary value problem for the subdomain:
\begin{align}\label{BVP_sub}
    \begin{cases}
                   \Div \mathbf{u}_\epsilon^s = f \qquad &\text{in $\Omega_\epsilon^s$},  \\
             \mathbf{u}_\epsilon^s +  \nabla p_\epsilon^s  = 0 \qquad &\text{in $\Omega_\epsilon^s \backslash\gamma$},   \\
		   t_f^{-1}\mathbf{u}_\epsilon^s \cdot \mathbf{n}  = \llbracket p_\epsilon^s \rrbracket  &\text{on $\gamma$}, \\
		            \mathbf{u}_\epsilon^s\cdot \mathbf{\nu} = \gamma_0(\mathbf{u}_\epsilon \cdot \mathbf{\nu}) &\text{on $\Gamma_{\epsilon}^{s,+} \cup \Gamma_{\epsilon}^{s,-}$} \\
                   p_\epsilon^s = \gamma_0p_\epsilon &\text{on $\partial \Omega_\epsilon^{s} \backslash (\Gamma_{\epsilon}^{s,+} \cup \Gamma_{\epsilon}^{s,-})$},
    \end{cases}
\end{align}
where $\gamma_0$ is the trace operator defined on $\Omega_\epsilon^s$. We can solve this problem by using the mixed variational formulation similar to the one in \eqref{vari_mix}. But now the solution space is changed to $\mathcal{H}_\epsilon = \tilde H(\DDiv, \Omega_\epsilon^s) \times H^1(\Omega_\epsilon^s \backslash \gamma)$. Then we can define our modified solution for given $\epsilon$ by
\begin{align}\label{modify}
    \tilde{p}_\epsilon = \begin{cases}
    p_\epsilon(\mathbf{x}) &\text{if $\mathbf{x}\in \Omega \backslash \Omega_\epsilon^s$},\\
    p_\epsilon^s(\mathbf{x}) &\text{if $\mathbf{x}\in\Omega_\epsilon^s$}, \qquad \qquad
    \end{cases}
    \tilde{\mathbf{u}}_\epsilon = \begin{cases}
    \mathbf{u}_\epsilon(\mathbf{x}) &\text{if $\mathbf{x}\in \Omega \backslash \Omega_\epsilon^s$},\\
    \mathbf{u}_\epsilon^s(\mathbf{x}) &\text{if $\mathbf{x}\in\Omega_\epsilon^s$}.
    \end{cases}
\end{align}

\subsection{Formulation of \texorpdfstring{$\delta_\epsilon$}{Lg}}
\begin{figure}
    \centering
    \begin{tikzpicture}[thick,scale=0.8, every node/.style={scale=0.8}]
    \pic at (1,0) {annotated cuboid={width=40, height=40, depth=40}};
    \draw[dashed] (-2.2,.8) -- (.8, 1.53);
    \draw[dashed] (-2.2,.8) -- (-1,-4);
    \draw[dashed] (-1,-4) -- (2, -3);
    \draw[dashed] (.8, 1.53) -- (2, -3);
       \draw[line width=1.5pt,blue,-stealth](-0.3,-1.5)--(0.,-1.35) node[anchor=south west]{\large$\mathbf{n}$};
    \node[] at (0, -2) {\Large$\gamma$};
    \node[] at (7, -2) {\Large$\gamma$};
    \draw[c3] (-1.3,0.3) -- (.4, .75);
    \draw[c3] (-1.3,0.3) -- (-.5,-3);
    \draw[c3] (-.5,-3) -- (1.2, -2.5);
    \draw[c3] (.4, .75) -- (1.2, -2.5);
    \fill[pattern=flexible north east lines, pattern density = 8pt, pattern line width=.5pt,pattern color=c3] (-1.3,0.3) -- (.4, .75) -- (1.2, -2.5) -- (-.5,-3) --  (-1.3,0.3) ;
    \draw [c3] (6,-2.3) rectangle (10,0.3);
    \fill[pattern=flexible north east lines, pattern density = 8pt, pattern line width=.5pt,pattern color=c3]  (6,-2.3) rectangle (10,0.3);
    \draw[line width=1.5pt,blue,-stealth](8,-1.5)--(8,-1);
    \draw[line width=1.5pt,blue,-stealth](8,-1.5)--(8.5,-1.5);
    \node[blue] at (8, -0.8) {\large$\mathbf{\tau}_2$};
    \node[blue] at (8.8, -1.5) {\large$\mathbf{\tau}_1$};
    \draw [black,dashed] (5,-3) rectangle (11,1);
    \node[] at (0, -4.5) {(a) three-dimensional view};
    \node[] at (8, -4.5) {(b) two-dimensional view};
    \end{tikzpicture}
    \caption{An illustration of the coordinates used for definition of $\delta_\epsilon$ in three-dimensional space. }
    \label{3d_example}
\end{figure}
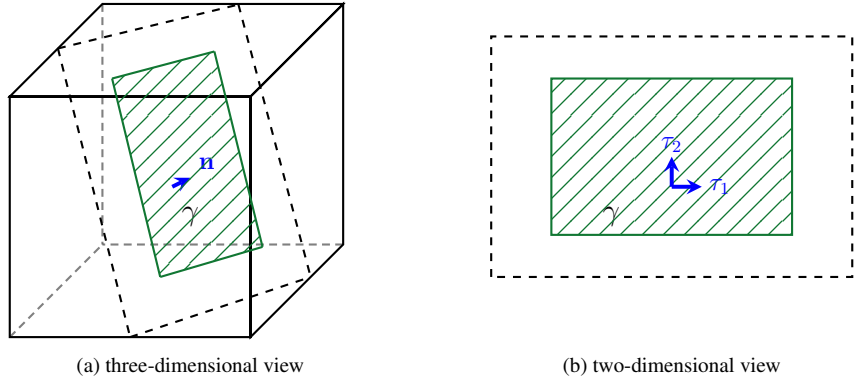

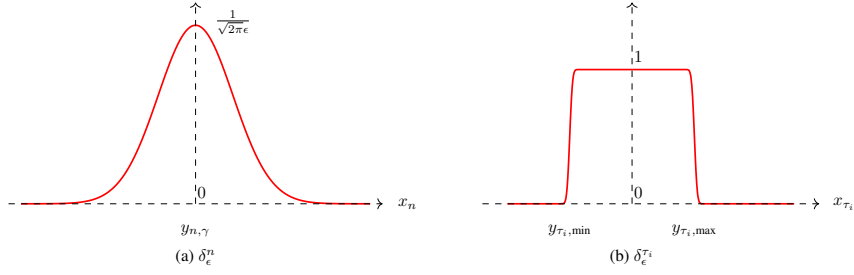
\begin{figure}
\centering
\begin{tikzpicture}[thick,scale=0.8, every node/.style={scale=0.8}]
\begin{axis}[
			axis line style={draw=none},
    		tick style={draw=none},
          xmin = -5, xmax=3,
          ymin = -.5,  ymax=1.75,
          xticklabels={,,},
          yticklabels={,,},
          width=1.1\textwidth, height=0.4\textwidth]
\addplot[domain=-4.4:-1.6, red,thick,smooth,samples=100] {1.33*e^(-((x+3.)/0.4243)^2.)};
\draw[->,dashed,black] (axis cs: -4.5,0)--(axis cs: -1.5, 0);
\draw[->,dashed,black] (axis cs: -3,0)--(axis cs: -3, 1.5);
\node at (axis cs:-2.7, 1.35) {$\frac{1}{\sqrt{2\pi}\epsilon}$};
\node at (axis cs:-2.95, 0.08) {$0$};
\node at (axis cs:-3, -0.2) {$y_{n,\gamma}$};
\node at (axis cs:-3, -0.4) {(a) $\delta_\epsilon^n$};
\addplot[domain=-0.5:.5, red,thick,smooth, samples=100] {.5*(1. + erf(40.*x))};
\addplot[domain=.5:1.8, red,thick,smooth, samples=100] {.5*(1. + erf(40.*(1.-x)))};
\draw[->,dashed,black] (axis cs: -.7,0)--(axis cs: 2, 0);
\draw[->,dashed,black] (axis cs: 0.5,0)--(axis cs: 0.5, 1.5);
\node at (axis cs:0.0, -0.2) {$y_{\tau_i,\text{min}}$};
\node at (axis cs:1.0, -0.2) {$y_{\tau_i,\text{max}}$};
\node at (axis cs:0.55, 1.1) {$1$};
\node at (axis cs:0.55, 0.08) {$0$};
\node at (axis cs:-1.3, 0) {$x_n$};
\node at (axis cs:2.2, 0) {$x_{\tau_i}$};
\node at (axis cs:.5, -0.4) {(b) $\delta_\epsilon^{\tau_i}$};
\end{axis}
\end{tikzpicture}
	\caption{An illustration of functions used for definition of $\delta_\epsilon$.}
	\label{delta}
\end{figure}

We deploy coordinates based on the fualt, i.e.,  $\{\mathbf{n}, \mathbf{\tau}_1,\mathbf{\tau}_2,\ldots,\mathbf{\tau}_{d-1} \}$. Assume that the fault is rectangle in three-dimensional space with $\tau_1$ and  $\tau_2$ defined along the length and width of the fault (see \autoref{3d_example}). Suppose that each point $\mathbf{y}$ on the fault could be expressed as a tuple $(y_{n,\gamma}, y_{\tau_1}, \ldots, y_{\tau_{d-1}})$, where $y_{n,\gamma}$ is a constant for all points on the fault. Denote the minimum and maximum values of the projection in $\tau_i$ direction by $y_{\tau_i,\text{min}}$ and $y_{\tau_i,\text{max}}$. We decompose $\delta_\epsilon$ into $d$ components in each direction of the coordinate and define it by
\begin{align*}
    \delta_\epsilon(\mathbf{x}) = \delta_\epsilon^{n}(\mathbf{x})\prod_{i=1}^{d-1} \delta_\epsilon^{\tau_i}(\mathbf{x}) \qquad \forall \mathbf{x}=(x_n,x_{\tau_i}, \ldots, x_{\tau_{d-1}}) \in \Omega,
\end{align*}
where
\begin{align*}
    \delta_\epsilon^n(\mathbf{x}) = \frac{1}{\sqrt{2\pi}\epsilon}\exp{\Big(-\frac{(x_n-y_{n,\gamma})^2}{2\epsilon^2}\Big)},
\end{align*}
and 
\begin{align*}
    \delta_\epsilon^{\tau_i} (\mathbf{x}) =\nolinebreak
   \frac{1}{4} \Big[1+\erf\Big(\frac{x_{\tau_i} - y_{\tau_i, \text{min}}}{\sqrt{2\pi}\epsilon_\tau}\Big) \Big] \nolinebreak
   \Big[1+\erf\Big(\frac{y_{\tau_i, \text{max}} - x_{\tau_i}}{\sqrt{2\pi}\epsilon_\tau}\Big) \Big].
\end{align*}
\autoref{delta} shows an example of the shape of these two kind of functions. Since both $\delta_n$ and $\delta_{\tau_i}(1\leq i \leq d-1)$ are in $ C^\infty$-functions, $\delta_\epsilon \in C^\infty(\Omega)$. Considering both $\epsilon$ and $\epsilon_\tau$ are sufficiently small opposed to the scale of $\Omega$, we could assume that $\delta_\epsilon$ has compact support in $\Omega$ in practice.

\subsection{Implementation on FEniCS}
For an efficient implementation of both the mixed method and our new method, we chose to use Gmsh \cite{gmsh} for mesh generation and base our code on FEniCS finite element platform \cite{Logg2012}. 

To solve the mixed variational problem \eqref{vari_mix} and the mixed problem associated with subdomains around faults \eqref{BVP_sub} in the new method, we use Discontinuous Galerkin basis functions ('DG' element in FEniCS) with degree of 0 to define the pressure function space, and use Raviart–Thomas basis functions ('RT' element in FEniCS) with degree of 1 to construct the velocity function space. The solution space is the product of these two spaces with degree of freedom (dof) equal to the sum of the number of elements and faces of the mesh.

To solve approximate pressure $p_\epsilon$ in the continuous variational problem \eqref{weak_2} associated with the whole domain of the new method, we use the Lagrange basis function ('Lagrange' element in FEniCS) with degree of 1 to define the pressure function space. Note that dof associated with this problem is the number of vertices of the mesh. For the corresponding velocity $\mathbf{u}_\epsilon$, we still use the Lagrange basis function with degree of 1 to construct the vector function space.

As for the linear solver, we use the iterative method GMRES. For the experiments showed in this paper, absolute tolerance \num{1.e-8} is set for GMRES. In addition, incomplete LU factorization with 1 level of fill, i.e. ILU(1), is used as preconditioner for the linear solver.

\section{Numerical Expriments} \label{results}
In this section, two numerical experiments are conducted to validate the new method. We first consider a two-dimensional domain where both high-transmissibility fault and low-transmissibility fault are tested repectively to examine the accuracy of the new method (\autoref{2D}). Then, a three-dimensional test is performed to show the efficiency of the new method (\autoref{3D}). For each test, ther performance of the new method is compared with the mixed method under different mesh resolutions.

For simplicity, we assume $f = 0$, $\frac{\rho}{\mu}=1$, and $\kappa$ as identity matrix in \eqref{BVP} in each test. The errors of each test is measured in the following way. The mixed method is performed at high resolution with fine grids, the solution of which is used to offer a reference as the ground truth solution. Then the $L^2$ errors of the pressure and velocity solutions are given by
\begin{align*}
    \Vert e_p^h \Vert_{L^2} = \Vert p^h - p^g \Vert_{L^2(\Omega)}, \qquad \nolinebreak
    \Vert \mathbf{e_u^h} \Vert_{L^2} = \Vert \mathbf{u^h} - \mathbf{u^g} \Vert_{L^2(\Omega)},
\end{align*}
where $p^h$ and $\mathbf{u}^h$ are the pressure and velocity results using the derived new method under mesh size $h$, $p^g$ and $\mathbf{u^g}$ are the pressure and velocity results of the ground truth.

\subsection{2D test}\label{2D}
\paragraph{Setup.} 

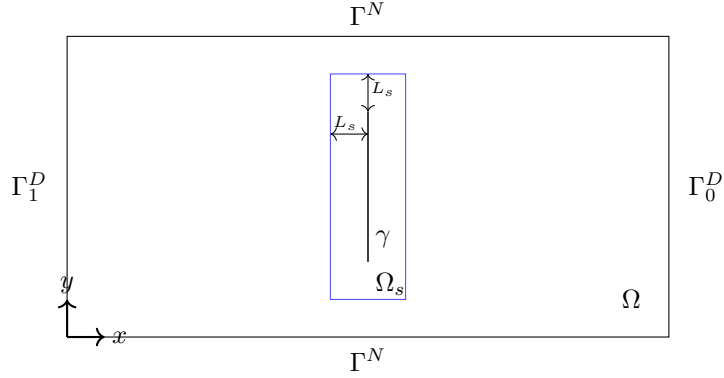
\begin{figure}[t!]
\begin{center}
  \begin{tikzpicture}
    \draw[black] (-4, -2) -- (4,-2) -- (4,2) -- (-4,2) -- (-4,-2) ;
    \draw[black, line width = 0.2mm] (0,-1) to (0,1);
    \draw[blue!70] (-0.5, -1.5) -- (0.5,-1.5) -- (.5,1.5) -- (-.5,1.5) -- (-.5,-1.5) ;
    \draw[<->] (-0.5, 0.7) to (0,0.7);
    \draw[<->] (-0, 1) to (0,1.5);
    \node[] at (-0.3, 0.85)   (a) {\tiny$L_s$};
     \node[] at (0.2, 1.3)   (a) {\tiny$L_s$};
    \node[] at (3.5,-1.5)   (a) {$\Omega$};
    \node[] at (4.5,0)   (b) {$\Gamma_0^D$};
    \node[] at (-4.5,0)   (b) {$\Gamma_1^D$};
    \node[] at (0., 2.3)   (b) {$\Gamma^N$};
    \node[] at (0., -2.3)   (b) {$\Gamma^N$};
    \node[] at (0.2,-0.7)   (d) {$\gamma$};
    \node[] at (0.3,-1.3)   (d) {$\Omega_s$};
    \draw[thick, ->] (-4, -2) to (-3.5,-2);
    \draw[thick, ->] (-4, -2) to (-4,-1.5);
    \node[] at (-3.3, -2)   (a) {$x$};
    \node[] at (-4, -1.3)   (a) {$y$};

  \end{tikzpicture}
  \end{center}
  \caption{Configuration of the 2D test.}
  \label{2D_setup}
\end{figure}

We consider a two-dimensional rectangular domain with width $L_x = 2.0$ and length $L_y = 1.0$ as showed in \autoref{2D_setup}. A fault $\gamma$, as a line with length of $0.4$, is located in the middle of the domain with normal direction in $x$. We split the boundary as three parts: $\Gamma_1^D = \{(0,y): 0 \leq y \leq L_y \}$, $\Gamma_0^D = \{ (1,y): 0 \leq y \leq L_y \}$, and $\Gamma^N = \{ (x,y): 0 \leq x \leq L_x, \text{y = 0 or 1} \}$. The boundary conditions are given by
\begin{align*}
\begin{cases}
    p = 1 \qquad \text{on $\Gamma_1^D$},\\
    p = 0  \qquad \text{on $\Gamma_0^D$},\\
    \nabla \mathbf{u}\cdot \mathbf{\nu} = 0 \qquad \text{on $\Gamma^N$}.
\end{cases}
\end{align*}

Meshes used for this test are set in the following way. We denote the mesh size on the domain boundary $\partial \Omega$ by $h$, the mesh size on the boundary of the subdomain $\partial \Omega_s$ by $h_s$, and mesh size on the fault $\gamma$ by $h_f$. The three mesh sizes have the relation of $h:h_s:h_f = 5:2:2$. Moreover, we define the distance between subdomain boundary to the fault as $L_s$ in \autoref{2D_setup} and set $L_s = 20h_f$.

Five different meshes are performed for the new method and the mixed method for the cases with $t_f=2.0$ and $t_f = 0.02$ respectively. For the new method, we found that $\epsilon=3h_f$ is the optimal when $t_f = 2.0$ and $\epsilon = 1h_f$ is the optimal when $t_f = 0.02$. We show results of $\epsilon=3h_f$, 4$h_f$ and $5h_f$ for the new method when $t_f = 2.0$ and $\epsilon=1h_f$, 2$h_f$ and $3h_f$ for the new method when $t_f = 0.02$. The ground truth solutions are from the mixed method with $h=2.5\mathrm{E}{-3}$.

\begin{figure}[t!]
\centering
\begin{tikzpicture}[every plot/.append style={thick}]
\pgfplotsset{
scale only axis,
every axis legend/.append style={
legend columns=4
}
}
\begin{groupplot}[
xlabel near ticks, ylabel near ticks,
width=4cm,
height=3.5cm,
group style={
group name=2d, group size=3 by 2,
horizontal sep=.5cm, vertical sep=1.2cm,
y descriptions at=edge left,
}
]
\nextgroupplot[legend cell align={left},
legend style={nodes={scale=1.2, transform shape},
anchor=south west,
draw = none,line width=1.0pt},
legend to name={CommonLegend},
xlabel=$h$,
xmin = 0.002, xmax = 0.2,
xmode = log,
ymax = 2.e-2,ymin = 5.e-4,
ylabel=$\Vert e_p^h \Vert_{L^2}$,
ymode=log,
log basis y={10},
xlabel near ticks, ylabel near ticks,
grid = both]
\addplot[red,mark= square, mark options={solid}] table [x=h ,y=ep]
{fig/2D_mix_error_2.txt};
\addplot[c3,mark= o,mark options={solid}] table [x=h ,y=ep_1]
{fig/2D_con_error_2.txt};
\addplot[blue!70,mark= asterisk, mark options={solid}] table [x=h ,y=ep_2]
{fig/2D_con_error_2.txt};
\addplot[black,mark= diamond,mark options={solid}] table [x=h ,y=ep_3]
{fig/2D_con_error_2.txt};
\draw[black,thick] (axis cs:0.004,0.004) -- node[left= 0.05cm] {$1$} (axis cs:0.004,0.008);
\draw[black,thick] (axis cs:0.004,0.008) -- node[above= 0.05cm] {$2$} (axis cs:0.016,0.008);
\draw[black,thick] (axis cs:0.016,0.008) --  (axis cs:0.004,0.004);

\draw[black,thick] (axis cs:0.03,0.001) -- node[below= 0.05cm] {$1$} (axis cs:0.06,0.001);
\draw[black,thick] (axis cs:0.06,0.001) -- node[right= 0.05cm] {$1$} (axis cs:0.06,0.002);
\draw[black,thick] (axis cs:0.06,0.002) --  (axis cs:0.03,0.001);

\addlegendimage{/pgfplots/refstyle=plot_one}
\addlegendentry{\scriptsize{mixed method}}
\addlegendentry{\scriptsize{$\epsilon = 3h_f$}}
\addlegendentry{\scriptsize{$\epsilon = 4h_f$}}
\addlegendentry{\scriptsize{$\epsilon = 5h_f$}}
\label{plot_one}

\nextgroupplot[ legend cell align={left},
xlabel near ticks, ylabel near ticks,
legend style={at={(axis cs:1000,0.015)},nodes={scale=0.8, transform shape},anchor=south west, draw = none,line width=1.0pt},
xlabel=dof,
xmin = 1000, xmax = 500000,
xmode = log,
ymax = 2.e-2,ymin = 5.e-4,
ymode=log,
log basis y={10},
grid = both]
\addplot[red,mark= square, mark options={solid}] table [x=dof ,y=ep]
{fig/2D_mix_error_2.txt};
\addplot[c3,mark= o,mark options={solid}] table [x=dof,y=ep_1]
{fig/2D_con_error_2.txt};
\addplot[blue!70,mark= asterisk, mark options={solid}] table [x=dof ,y=ep_2]
{fig/2D_con_error_2.txt};
\addplot[black,mark= diamond,mark options={solid}] table [x=dof ,y=ep_3]
{fig/2D_con_error_2.txt};
\nextgroupplot[ legend cell align={left},
xlabel near ticks, ylabel near ticks,
legend style={at={(axis cs:0.01,0.015)},nodes={scale=0.8, transform shape},anchor=south west, draw = none,line width=1.0pt},
xlabel=time,
xmin = 0.01, xmax = 100,
xmode = log,
ymax = 2.e-2,ymin = 5.e-4,
ymode=log,
log basis y={10},
grid = both]
\addplot[red,mark= square, mark options={solid}] table [x=t ,y=ep]
{fig/2D_mix_error_2.txt};
\addplot[c3,mark= o,mark options={solid}] table [x=t_1,y=ep_1]
{fig/2D_con_error_2.txt};
\addplot[blue!70,mark= asterisk, mark options={solid}] table [x=t_2 ,y=ep_2]
{fig/2D_con_error_2.txt};
\addplot[black,mark= diamond,mark options={solid}] table [x=t_3 ,y=ep_3]
{fig/2D_con_error_2.txt};

\nextgroupplot[ legend cell align={left},
legend style={at={(axis cs:0.002,0.04)},nodes={scale=0.8, transform shape},anchor=south west, draw = none,line width=1.0pt},
xlabel=$h$,
xmin = 0.002, xmax = 0.2,
xmode = log,
ymax = 1.e-1,ymin = 4.e-3,
ylabel=$\Vert \mathbf{e_u^h} \Vert_{L^2}$,
ymode=log,
log basis y={10},
xlabel near ticks, ylabel near ticks,
grid = both]
\addplot[red,mark= square, mark options={solid}] table [x=h ,y=eu]
{fig/2D_mix_error_2.txt};
\addplot[c3,mark= o,mark options={solid}] table [x=h ,y=eu_1]
{fig/2D_con_error_2.txt};
\addplot[blue!70,mark= asterisk, mark options={solid}] table [x=h ,y=eu_2]
{fig/2D_con_error_2.txt};
\addplot[black,mark= diamond,mark options={solid}] table [x=h ,y=eu_3]
{fig/2D_con_error_2.txt};
\draw[black,thick] (axis cs:0.03,0.008) -- node[below= 0.03cm] {$2$} (axis cs:0.12, 0.008);
\draw[black,thick] (axis cs:0.12,0.008) -- node[right= 0.05cm] {$1$} (axis cs:0.12,0.016);
\draw[black,thick] (axis cs:0.12,0.016) --  (axis cs:0.03,0.008);

\nextgroupplot[ legend cell align={left},
xlabel near ticks, ylabel near ticks,
legend style={at={(axis cs:1000,0.04)},nodes={scale=0.8, transform shape},anchor=south west, draw = none,line width=1.0pt},
xlabel=dof,
xmin = 1000, xmax = 500000,
xmode = log,
ymax = 1.e-1,ymin = 4.e-3,
ymode=log,
log basis y={10},
grid = both]
\addplot[red,mark= square, mark options={solid}] table [x=dof ,y=eu]
{fig/2D_mix_error_2.txt};
\addplot[c3,mark= o,mark options={solid}] table [x=dof,y=eu_1]
{fig/2D_con_error_2.txt};
\addplot[blue!70,mark= asterisk, mark options={solid}] table [x=dof ,y=eu_2]
{fig/2D_con_error_2.txt};
\addplot[black,mark= diamond,mark options={solid}] table [x=dof ,y=eu_3]
{fig/2D_con_error_2.txt};

\nextgroupplot[ legend cell align={left},
xlabel near ticks, ylabel near ticks,
legend style={at={(axis cs:0.01,0.04)},nodes={scale=0.8, transform shape},anchor=south west, draw = none,line width=1.0pt},
xlabel=time,
xmin = 0.01, xmax = 100,
xmode = log,
ymax = 1.e-1,ymin = 4.e-3,
ymode=log,
log basis y={10},
grid = both]
\addplot[red,mark= square, mark options={solid}] table [x=t ,y=eu]
{fig/2D_mix_error_2.txt};
\addplot[c3,mark= o,mark options={solid}] table [x=t_1,y=eu_1]
{fig/2D_con_error_2.txt};
\addplot[blue!70,mark= asterisk, mark options={solid}] table [x=t_2 ,y=eu_2]
{fig/2D_con_error_2.txt};
\addplot[black,mark= diamond,mark options={solid}] table [x=t_3 ,y=eu_3]
{fig/2D_con_error_2.txt};

\end{groupplot}
\path (2d c2r1.north west) -- node[above=0.1cm,font=\fontsize{12pt}{5pt}\selectfont]{\ref{CommonLegend}} (2d c2r1.north east);
\end{tikzpicture}
\caption{$L^2$ errors of pressure and velocity for the new method and the mixed method in the 2D test with $t_f = 2.0$.}
\label{2d_convergence_2}
\end{figure}
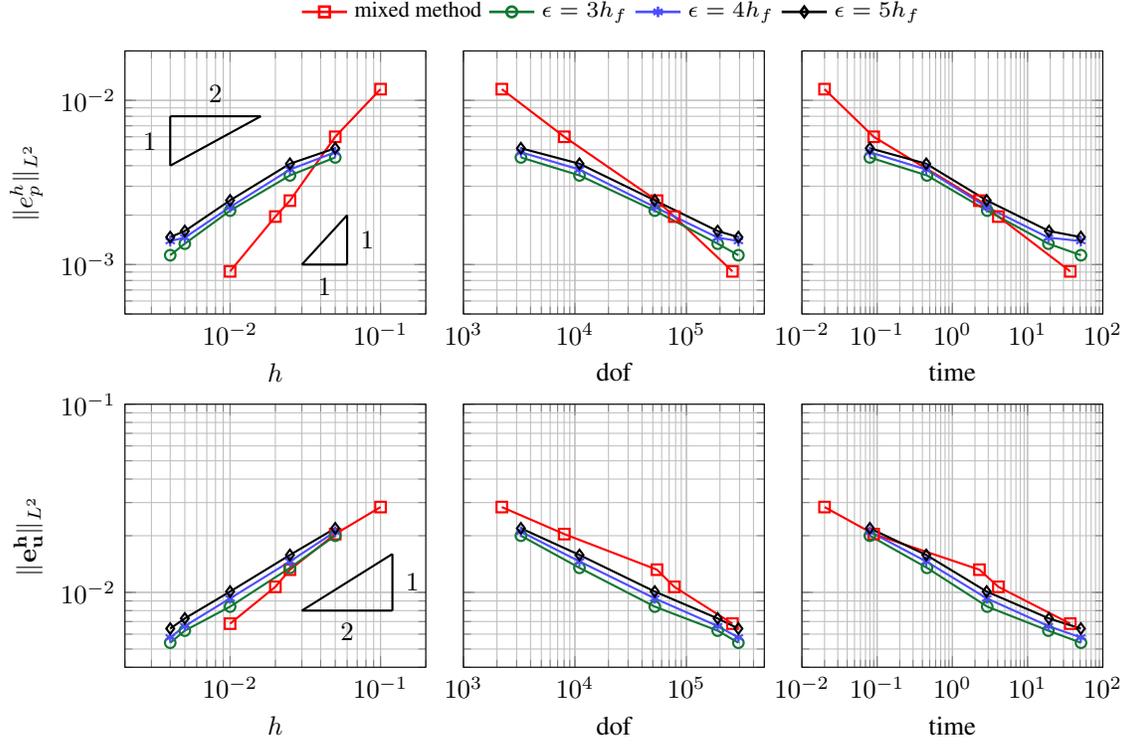

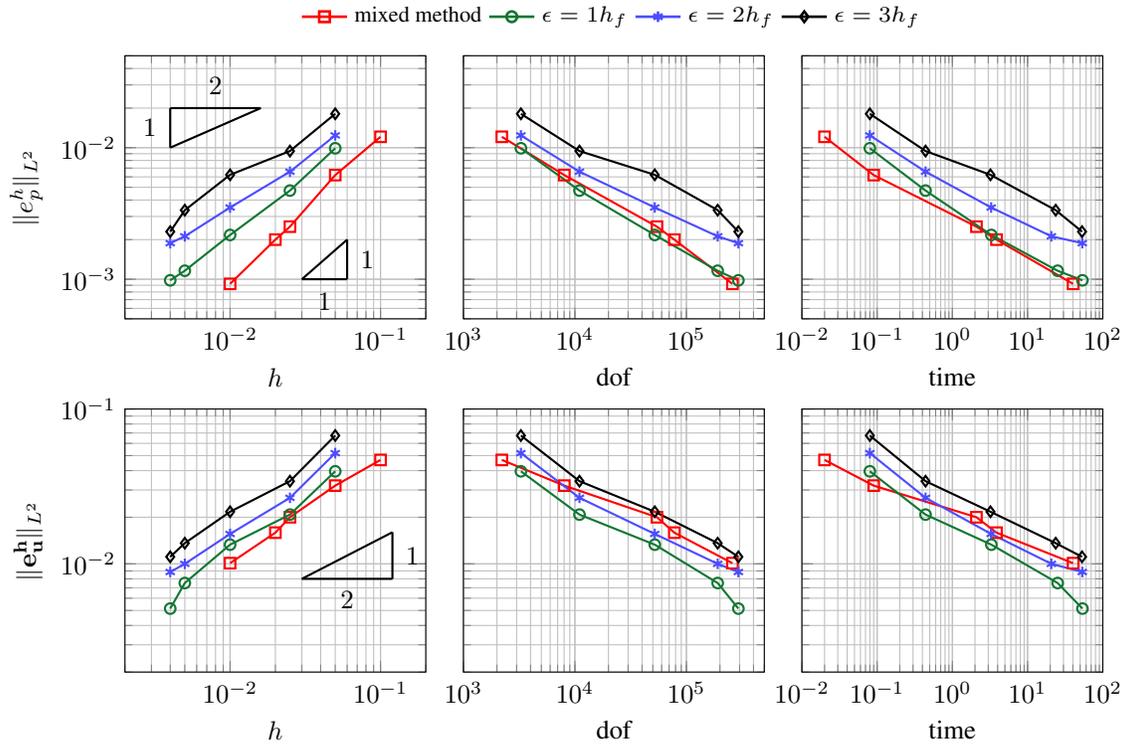
\begin{figure}[t!]
\centering
\begin{tikzpicture}[every plot/.append style={thick}]
\pgfplotsset{
scale only axis,
every axis legend/.append style={
legend columns=4
}
}
\begin{groupplot}[
xlabel near ticks, ylabel near ticks,
width=4cm,
height=3.5cm,
group style={
group name=2d, group size=3 by 2,
horizontal sep=.5cm, vertical sep=1.2cm,
y descriptions at=edge left,
}
]
\nextgroupplot[legend cell align={left},
legend style={nodes={scale=1.2, transform shape},
anchor=south west,
draw = none,line width=1.0pt},
legend to name={CommonLegend},
xlabel=$h$,
xmin = 0.002, xmax = 0.2,
xmode = log,
ymax = 5.e-2,ymin = 5.e-4,
ylabel=$\Vert e_p^h \Vert_{L^2}$,
ymode=log,
log basis y={10},
xlabel near ticks, ylabel near ticks,
grid = both]
\addplot[red,mark= square, mark options={solid}] table [x=h ,y=ep]
{fig/2D_mix_error_002.txt};
\addplot[c3,mark= o,mark options={solid}] table [x=h ,y=ep_1]
{fig/2D_con_error_002.txt};
\addplot[blue!70,mark= asterisk, mark options={solid}] table [x=h ,y=ep_2]
{fig/2D_con_error_002.txt};
\addplot[black,mark= diamond,mark options={solid}] table [x=h ,y=ep_3]
{fig/2D_con_error_002.txt};
\draw[black,thick] (axis cs:0.004,0.01) -- node[left= 0.05cm] {$1$} (axis cs:0.004,0.02);
\draw[black,thick] (axis cs:0.004,0.02) -- node[above= 0.05cm] {$2$} (axis cs:0.016,0.02);
\draw[black,thick] (axis cs:0.016,0.02) --  (axis cs:0.004,0.01);

\draw[black,thick] (axis cs:0.03,0.001) -- node[below= 0.05cm] {$1$} (axis cs:0.06,0.001);
\draw[black,thick] (axis cs:0.06,0.001) -- node[right= 0.05cm] {$1$} (axis cs:0.06,0.002);
\draw[black,thick] (axis cs:0.06,0.002) --  (axis cs:0.03,0.001);

\addlegendimage{/pgfplots/refstyle=plot_one}
\addlegendentry{\scriptsize{mixed method}}
\addlegendentry{\scriptsize{$\epsilon = 1h_f$}}
\addlegendentry{\scriptsize{$\epsilon = 2h_f$}}
\addlegendentry{\scriptsize{$\epsilon = 3h_f$}}
\label{plot_one}

\nextgroupplot[ legend cell align={left},
xlabel near ticks, ylabel near ticks,
legend style={at={(axis cs:1000,0.015)},nodes={scale=0.8, transform shape},anchor=south west, draw = none,line width=1.0pt},
xlabel=dof,
xmin = 1000, xmax = 500000,
xmode = log,
ymax = 5.e-2,ymin = 5.e-4,
ymode=log,
log basis y={10},
grid = both]
\addplot[red,mark= square, mark options={solid}] table [x=dof ,y=ep]
{fig/2D_mix_error_002.txt};
\addplot[c3,mark= o,mark options={solid}] table [x=dof,y=ep_1]
{fig/2D_con_error_002.txt};
\addplot[blue!70,mark= asterisk, mark options={solid}] table [x=dof ,y=ep_2]
{fig/2D_con_error_002.txt};
\addplot[black,mark= diamond,mark options={solid}] table [x=dof ,y=ep_3]
{fig/2D_con_error_002.txt};
\nextgroupplot[ legend cell align={left},
xlabel near ticks, ylabel near ticks,
legend style={at={(axis cs:0.01,0.015)},nodes={scale=0.8, transform shape},anchor=south west, draw = none,line width=1.0pt},
xlabel=time,
xmin = 0.01, xmax = 100,
xmode = log,
ymax = 5.e-2,ymin = 5.e-4,
ymode=log,
log basis y={10},
grid = both]
\addplot[red,mark= square, mark options={solid}] table [x=t ,y=ep]
{fig/2D_mix_error_002.txt};
\addplot[c3,mark= o,mark options={solid}] table [x=t_1,y=ep_1]
{fig/2D_con_error_002.txt};
\addplot[blue!70,mark= asterisk, mark options={solid}] table [x=t_2 ,y=ep_2]
{fig/2D_con_error_002.txt};
\addplot[black,mark= diamond,mark options={solid}] table [x=t_3 ,y=ep_3]
{fig/2D_con_error_002.txt};

\nextgroupplot[ legend cell align={left},
legend style={at={(axis cs:0.002,0.04)},nodes={scale=0.8, transform shape},anchor=south west, draw = none,line width=1.0pt},
xlabel=$h$,
xmin = 0.002, xmax = 0.2,
xmode = log,
ymax = 1.e-1,ymin = 2.e-3,
ylabel=$\Vert \mathbf{e_u^h} \Vert_{L^2}$,
ymode=log,
log basis y={10},
xlabel near ticks, ylabel near ticks,
grid = both]
\addplot[red,mark= square, mark options={solid}] table [x=h ,y=eu]
{fig/2D_mix_error_002.txt};
\addplot[c3,mark= o,mark options={solid}] table [x=h ,y=eu_1]
{fig/2D_con_error_002.txt};
\addplot[blue!70,mark= asterisk, mark options={solid}] table [x=h ,y=eu_2]
{fig/2D_con_error_002.txt};
\addplot[black,mark= diamond,mark options={solid}] table [x=h ,y=eu_3]
{fig/2D_con_error_002.txt};
\draw[black,thick] (axis cs:0.03,0.008) -- node[below= 0.03cm] {$2$} (axis cs:0.12, 0.008);
\draw[black,thick] (axis cs:0.12,0.008) -- node[right= 0.05cm] {$1$} (axis cs:0.12,0.016);
\draw[black,thick] (axis cs:0.12,0.016) --  (axis cs:0.03,0.008);

\nextgroupplot[ legend cell align={left},
xlabel near ticks, ylabel near ticks,
legend style={at={(axis cs:1000,0.04)},nodes={scale=0.8, transform shape},anchor=south west, draw = none,line width=1.0pt},
xlabel=dof,
xmin = 1000, xmax = 500000,
xmode = log,
ymax = 1.e-1,ymin = 2.e-3,
ymode=log,
log basis y={10},
grid = both]
\addplot[red,mark= square, mark options={solid}] table [x=dof ,y=eu]
{fig/2D_mix_error_002.txt};
\addplot[c3,mark= o,mark options={solid}] table [x=dof,y=eu_1]
{fig/2D_con_error_002.txt};
\addplot[blue!70,mark= asterisk, mark options={solid}] table [x=dof ,y=eu_2]
{fig/2D_con_error_002.txt};
\addplot[black,mark= diamond,mark options={solid}] table [x=dof ,y=eu_3]
{fig/2D_con_error_002.txt};

\nextgroupplot[ legend cell align={left},
xlabel near ticks, ylabel near ticks,
legend style={at={(axis cs:0.01,0.04)},nodes={scale=0.8, transform shape},anchor=south west, draw = none,line width=1.0pt},
xlabel=time,
xmin = 0.01, xmax = 100,
xmode = log,
ymax = 1.e-1,ymin = 2.e-3,
ymode=log,
log basis y={10},
grid = both]
\addplot[red,mark= square, mark options={solid}] table [x=t ,y=eu]
{fig/2D_mix_error_002.txt};
\addplot[c3,mark= o,mark options={solid}] table [x=t_1,y=eu_1]
{fig/2D_con_error_002.txt};
\addplot[blue!70,mark= asterisk, mark options={solid}] table [x=t_2 ,y=eu_2]
{fig/2D_con_error_002.txt};
\addplot[black,mark= diamond,mark options={solid}] table [x=t_3 ,y=eu_3]
{fig/2D_con_error_002.txt};

\end{groupplot}
\path (2d c2r1.north west) -- node[above=0.1cm,font=\fontsize{12pt}{5pt}\selectfont]{\ref{CommonLegend}} (2d c2r1.north east);
\end{tikzpicture}
\caption{$L^2$ errors of pressure and velocity for the new method and the mixed method in the 2D test with $t_f = 0.02$.}
\label{2d_convergence_002}
\end{figure}

\paragraph{Results.} In \autoref{2d_table_mix}, we investigate the simulation results for the mixed method under different mesh resolutions. Similarly, \autoref{2d_table_con_2} and \autoref{2d_table_con_002} report the results of the new method with $t_f = 2.0$ and $t_f=0.02$ respectively. CPU time for the mixed method represents the time used for solving the linear system. The CPU time is composed by the time used for solving the linear system of the pressure solution over the whole domain and the time for solving the mixed problem associated with the subdomain around the fault. We plot $L^2$ errors of pressure and velocity as functions of $h$, dof and time in \autoref{2d_convergence_2} and \autoref{2d_convergence_002}, respectively. We can see that both the two methods converge to the ground truth as the mesh size gets smaller. The pressure solution converges as $\mathcal{O}(h)$ in the mixed method, and seems to converge as $\mathcal{O}(h^{\frac{1}{2}})$ in our new method when $t_f=2.0$ and converge as $\mathcal{O}(h)$ when $t_f = 0.02$. As for the velocity solution, it converges as $\mathcal{O}(h^{\frac{1}{2}})$ for both of the two methods and both values of transmissibility. As for efficiency when $t_f = 2.0$, the new method could be faster than the mixed method at high pressure error (for example, higher than $3.0\mathrm{E}{-3}$) but could be slower at low pressure error; the new method seems to be a little bit faster than the mixed one from the plot of velocity error in \autoref{2d_convergence_2}. When $t_f=0.02$, the new method with $\epsilon=2h_f$ and $\epsilon=3h_f$ are both slower than the mixed method at given pressure error or velocity error. But the new method with $\epsilon=1h_f$ is very close to the mixed method at given pressure error, and could be faster than the mixed method at low velocity error (for example, lower than $3.0\mathrm{E}{-2}$ in \autoref{2d_convergence_002}.

\setlength{\extrarowheight}{5pt}
\begin{table}[!t]
\scriptsize
\centering
\caption{Summary of simulation results for the mixed method in the 2D test.}
\label{2d_table_mix}
\begin{tabular*}{\textwidth}{@{\extracolsep{3pt}}
P{0.12\textwidth-2\tabcolsep}
P{0.12\textwidth-2\tabcolsep}
P{0.12\textwidth-2\tabcolsep}
P{0.12\textwidth-2\tabcolsep}
P{0.12\textwidth-2\tabcolsep}
P{0.12\textwidth-2\tabcolsep}
P{0.12\textwidth-2\tabcolsep}
P{0.12\textwidth-2\tabcolsep}}
\toprule %
\multirow{2}{*}{h} &\multirow{2}{*}{dof} &\multicolumn{3}{c}{$t_f=2.0$} &\multicolumn{3}{c}{$t_f=0.02$}\\
\cline{3-5} \cline{6-8}
& &time &$\Vert e_p^h\Vert_{L^2}$ &$\Vert \mathbf{e_{\mathbf{u}}^h} \Vert_{L^2}$ &time &$\Vert e_p^h\Vert_{L^2}$ &$\Vert \mathbf{e_{\mathbf{u}}^h} \Vert_{L^2}$\\
\toprule
1.0E-01   &2210   &0.02   &1.17E-02   &2.84E-02   &0.02   &1.21E-02   &4.69E-02  \\
5.0E-02   &8035   &0.09   &6.00E-03   &2.04E-02   &0.09   &6.19E-03   &3.20E-02   \\
2.5E-02   &54455   &2.29   &2.45E-03   &1.32E-02   &2.09   &2.51E-03   &2.00E-02   \\
2.0E-02   &77965   &4.09   &1.96E-03   &1.07E-02   &3.86   &2.00E-03   &1.59E-02   \\
1.0E-02   &258970   &36.93   &9.09E-04   &6.82E-03   &40.27   &9.22E-04   &1.01E-02  \\
\bottomrule
\end{tabular*}
\end{table}

\begin{table}[!t]
\scriptsize
\centering
\vspace*{.5 cm}
\caption{Summary of simulation results of the new method in the 2D test with $t_f = 2.0$. The first number in "time" column represents the CPU time used for solving the approximate pressure and velocity of the whole domain, the second number is the time for the subdomain problem.}
\label{2d_table_con_2}
\begin{tabular*}{\textwidth}{@{\extracolsep{3pt}}
P{0.09\textwidth-2\tabcolsep}
P{0.08\textwidth-2\tabcolsep}
P{0.086\textwidth-2\tabcolsep}
P{0.087\textwidth-2\tabcolsep}
P{0.087\textwidth-2\tabcolsep}
P{0.086\textwidth-2\tabcolsep}
P{0.087\textwidth-2\tabcolsep}
P{0.087\textwidth-2\tabcolsep}
P{0.086\textwidth-2\tabcolsep}
P{0.087\textwidth-2\tabcolsep}
P{0.087\textwidth-2\tabcolsep}}
\toprule %
\multirow{2}{*}{h} &\multirow{2}{*}{dof} &\multicolumn{3}{c}{$\epsilon = 3h_f$} &\multicolumn{3}{c}{$\epsilon = 4h_f$} &\multicolumn{3}{c}{$\epsilon = 5h_f$}\\
\cline{3-5} \cline{6-8} \cline{9-11}
& &time &$\Vert e_p^h\Vert_{L^2}$ &$\Vert \mathbf{e_{\mathbf{u}}^h} \Vert_{L^2}$ &time &$\Vert e_p^h\Vert_{L^2}$ &$\Vert \mathbf{e_{\mathbf{u}}^h} \Vert_{L^2}$ &time &$\Vert e_p^h\Vert_{L^2}$ &$\Vert \mathbf{e_{\mathbf{u}}^h} \Vert_{L^2}$\\
\toprule
5.0E-02   &3279   &0.02+0.06   &4.5E-03   &2.0E-02   &0.02+0.06   &4.8E-03   &2.1E-02   &0.02+0.06   &5.1E-03   &2.2E-02   \\
2.5E-02   &10988   &0.10+0.36   &3.5E-03   &1.3E-02   &0.09+0.36   &3.8E-03   &1.5E-02   &0.10+0.35   &4.1E-03   &1.6E-02   \\
1.0E-02   &52035   &1.64+1.27   &2.1E-03   &8.4E-03   &1.60+1.27   &2.2E-03   &9.3E-03   &1.60+1.27   &2.5E-03   &1.0E-02   \\
5.0E-03   &190826   &13.94+4.99   &1.3E-03   &6.3E-03   &13.95+4.96   &1.5E-03   &6.6E-03   &14.29+4.94   &1.6E-03   &7.3E-03 \\
4.0E-03   &292641   &43.37+7.89   &1.1E-03   &5.4E-03   &43.11+7.93   &1.4E-03   &5.8E-03   &43.22+7.91   &1.5E-03   &6.4E-03   \\
\bottomrule
\end{tabular*}
\end{table}

\begin{table}[!t]
\scriptsize
\centering
\vspace*{.5 cm}
\caption{Summary of simulation results of the new method in the 2D test with $t_f=0.02$. The first number in "time" column represents the CPU time used for solving the approximate pressure and velocity of the whole domain, the second number is the time for the subdomain problem.}
\label{2d_table_con_002}
\begin{tabular*}{\textwidth}{@{\extracolsep{3pt}}
P{0.09\textwidth-2\tabcolsep}
P{0.08\textwidth-2\tabcolsep}
P{0.086\textwidth-2\tabcolsep}
P{0.087\textwidth-2\tabcolsep}
P{0.087\textwidth-2\tabcolsep}
P{0.086\textwidth-2\tabcolsep}
P{0.087\textwidth-2\tabcolsep}
P{0.087\textwidth-2\tabcolsep}
P{0.086\textwidth-2\tabcolsep}
P{0.087\textwidth-2\tabcolsep}
P{0.087\textwidth-2\tabcolsep}}
\toprule %
\multirow{2}{*}{h} &\multirow{2}{*}{dof} &\multicolumn{3}{c}{$\epsilon = 1h_f$} &\multicolumn{3}{c}{$\epsilon = 2h_f$} &\multicolumn{3}{c}{$\epsilon = 3h_f$}\\
\cline{3-5} \cline{6-8} \cline{9-11}
& &time &$\Vert e_p^h\Vert_{L^2}$ &$\Vert \mathbf{e_{\mathbf{u}}^h} \Vert_{L^2}$ &time &$\Vert e_p^h\Vert_{L^2}$ &$\Vert \mathbf{e_{\mathbf{u}}^h} \Vert_{L^2}$ &time &$\Vert e_p^h\Vert_{L^2}$ &$\Vert \mathbf{e_{\mathbf{u}}^h} \Vert_{L^2}$\\
\toprule
5.0E-02   &3279   &0.02+0.06   &9.9E-03   &4.0E-02   &0.02+0.06   &1.2E-02   &5.2E-02   &0.02+0.06   &1.8E-02   &6.7E-02   \\
2.5E-02   &10988   &0.10+0.33   &4.7E-03   &2.1E-02   &0.10+0.34   &6.6E-03   &2.7E-02   &0.11+0.33   &9.4E-03   &3.4E-02   \\
1.0E-02   &52035   &1.90+1.40   &2.2E-03   &1.3E-02   &1.88+1.40   &3.5E-03   &1.6E-02   &1.81+1.41   &6.2E-03   &2.2E-02   \\
5.0E-03   &190826   &19.64+5.48   &1.2E-03   &7.5E-03   &15.22+5.48   &2.1E-03   &1.0E-02   &18.19+5.46   &3.4E-03   &1.4E-02 \\
4.0E-03   &292641   &46.56+7.10   &9.8E-04   &5.1E-03   &46.59+6.77   &1.9E-03   &8.9E-03   &45.97+7.04   &2.3E-03   &1.1E-02   \\
\bottomrule
\end{tabular*}
\end{table}

\begin{figure}[t!]
\centering
\begin{tikzpicture}
\node[inner sep=0pt] (slice2) at (7, 0)
{\includegraphics[width=6cm,height = 4cm]{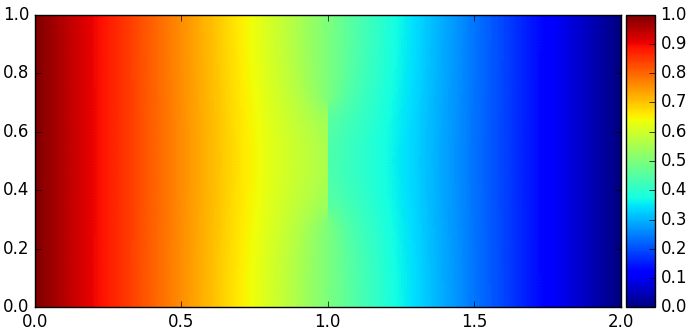}};
\node[inner sep=0pt] (slice1) at (0, 0)
{\includegraphics[width=6cm,height = 4cm]{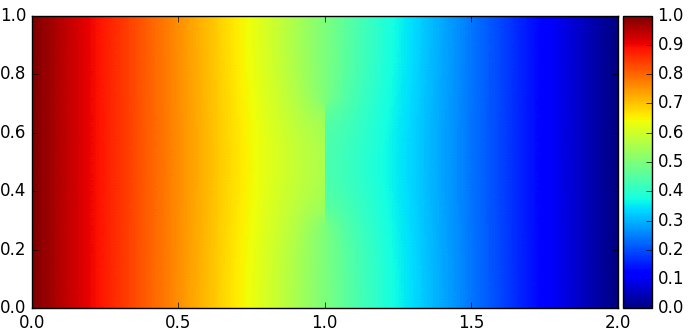}};
\node[inner sep=0pt] (slice3) at (7, -4.5)
{\includegraphics[width=6cm,height = 4cm]{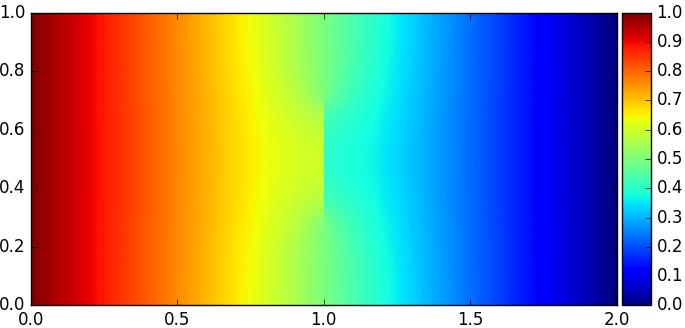}};
\node[inner sep=0pt] (slice4) at (0, -4.5)
{\includegraphics[width=6cm,height = 4cm]{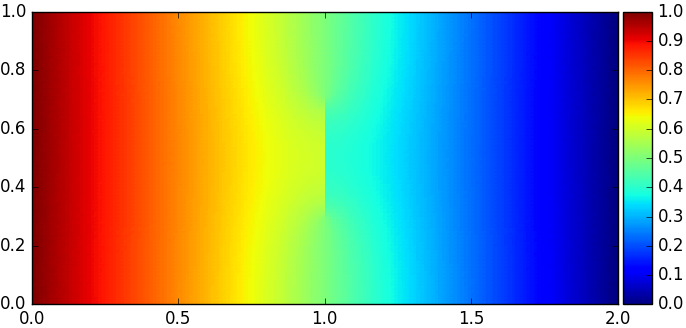}};
\end{tikzpicture}
\caption{Pressure solutions from the ground truth (left) and the new method of mesh with $h =5.0\mathrm{E}{-3}$ and $\epsilon = 3h_f$ (right) in the 2D test when $t_f = 2.0$ (upper) and $t_f=0.02$ (lower).}
\label{2d_p}
\end{figure}

\begin{figure}[!t]
\centering
\begin{tikzpicture}
\node[inner sep=0pt] (slice2) at (7, 0)
{\includegraphics[width=6cm,height = 4cm]{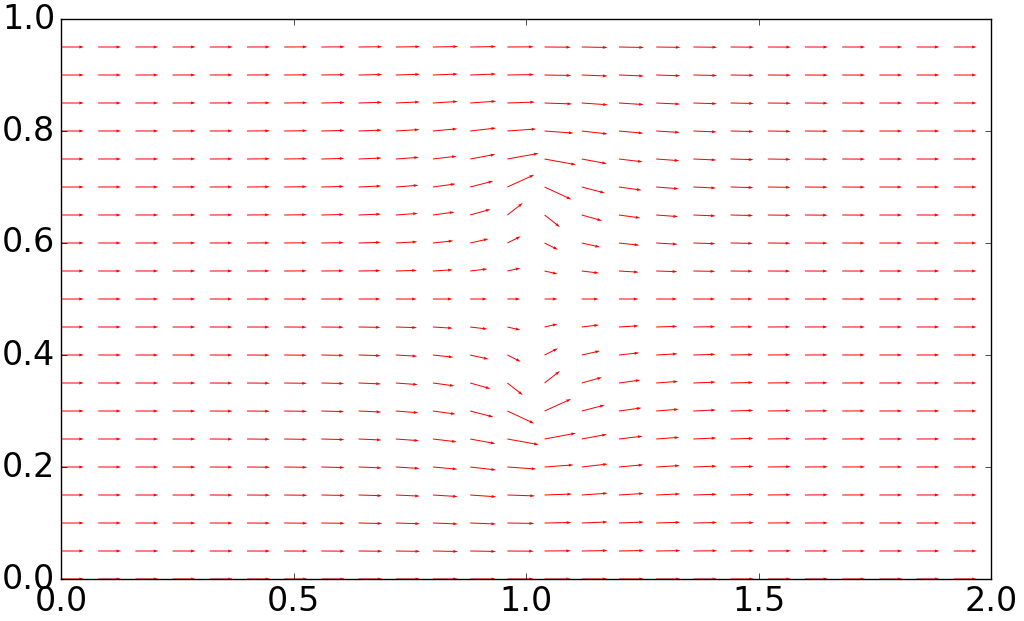}};
\node[inner sep=0pt] (slice1) at (0, 0)
{\includegraphics[width=6cm,height = 4cm]{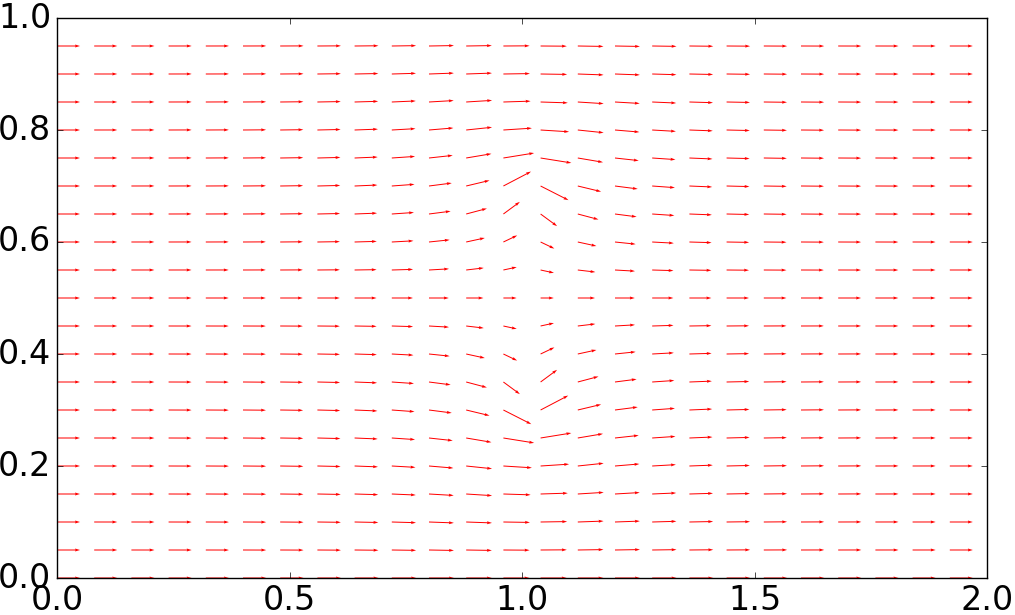}};
\node[inner sep=0pt] (slice3) at (7, -4.5)
{\includegraphics[width=6cm,height = 4cm]{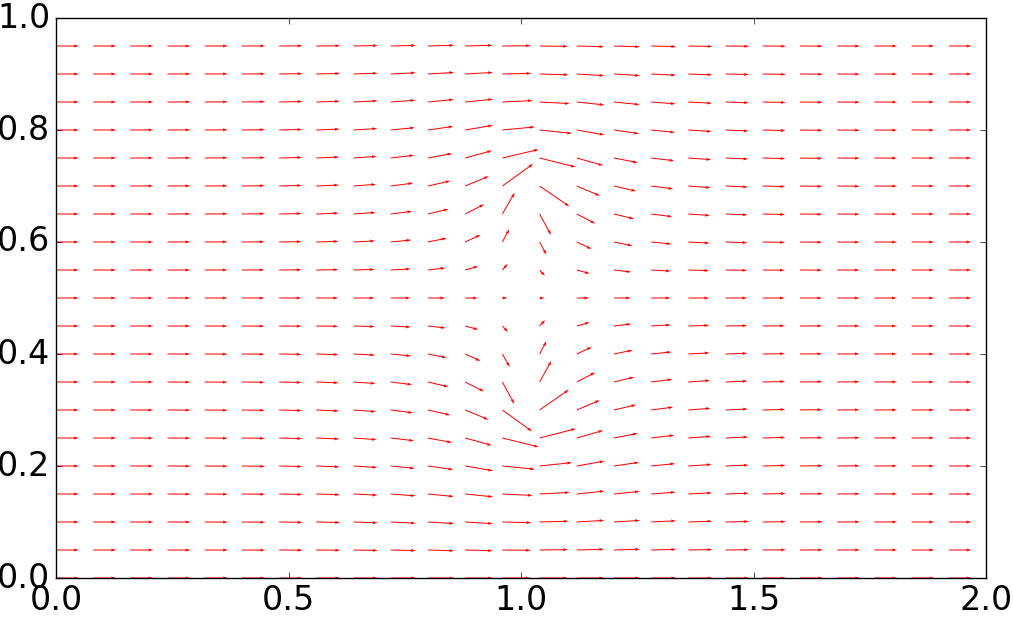}};
\node[inner sep=0pt] (slice4) at (0, -4.5)
{\includegraphics[width=6cm,height = 4cm]{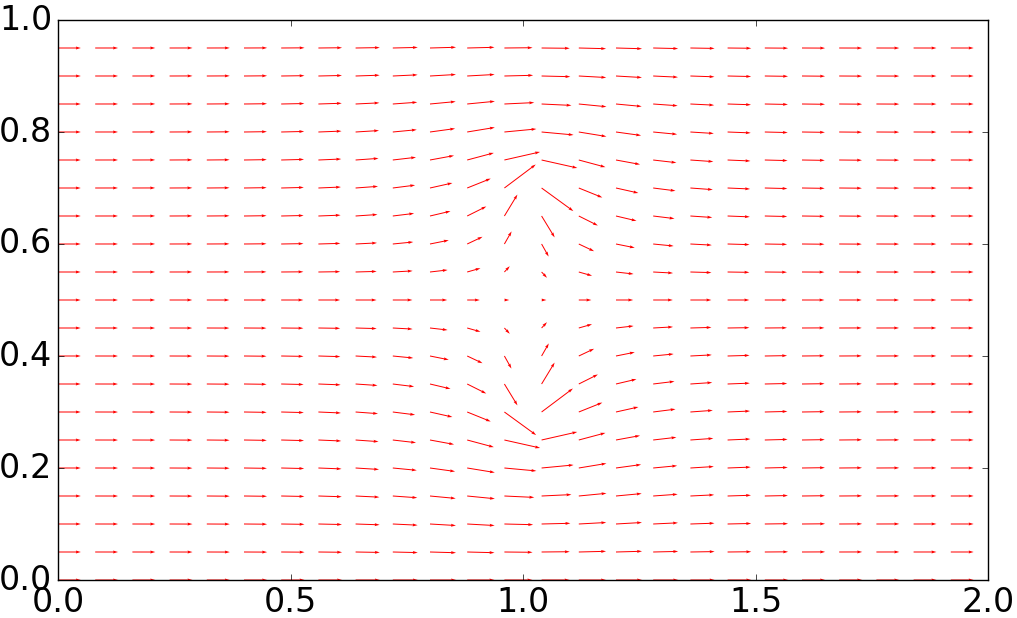}};
\end{tikzpicture}
\caption{Velocity solutions from the ground truth (left) and the new method of mesh with $h =5.0\mathrm{E}{-3}$ and $\epsilon = 3h_f$ (right) in the 2D test when $t_f = 2.0$ (upper) and $t_f=0.02$ (lower).}
\label{2d_u}
\end{figure}

Pressure and velocity solutions of the ground truth and one test of the mixed method are shown in \autoref{2d_p} and \autoref{2d_u}. In \autoref{2d_plot}, we plotted the pressure and the velocity in the normal direction along the line $\{(x,0.5): 0\leq x \leq L_x\}$. One can observe that both of the pressure and velocity solutions of the new method match the ground truth solutions in the whole domainm, which validates the accuracy of the new method. Moreover, lower transmissibility causes higher pressure jump and lower normal component of velocity along the fault. For the case of $t_f = 0.02$, the fault acts as a barrier for the flow such that nearly no flow could transport across the fault.

\begin{figure}[t!]
	\centering
	\begin{tikzpicture}[every plot/.append style={thick}]
	\begin{groupplot}[
	    footnotesize,
	    xlabel near ticks, ylabel near ticks,
        width=7cm,
        height=6cm,
        group style={
                     group name=case1, group size=2 by 1,
                     horizontal sep=2cm, vertical sep=2cm, 
                     },
        xlabel=i\textsuperscript{th} Largest Eigenvalue,
        ylabel=Eigenvalue,
      ]
     \nextgroupplot[ legend cell align={left},
   						legend style={at={(axis cs:32,300)},anchor=south west, draw = none,line width=1.0pt},
                        xmin = 0, xmax = 80,
                        xtick ={0,20,...,80},
                     xlabel near ticks, ylabel near ticks,
                        ymode = log,
                        ymin = 0, ymax = 100000]

                        \addplot[red] table [x=x ,y=tf1] {plot/mixed.txt};
                        \addplot[c3] table [x=x ,y=tf2] {plot/mixed.txt};
                        \addplot[blue!70] table [x=x ,y=tf3] {plot/mixed.txt};
                        \addplot[black] table [x=x ,y=tf4] {plot/mixed.txt};
                        
            \addlegendimage{/pgfplots/refstyle=plot_two}
            \addlegendentry{$t_f=2.0\mathrm{E}{0}$}
            \addlegendentry{$t_f=2.0\mathrm{E}{-1}$}
            \addlegendentry{$t_f=2.0\mathrm{E}{-2}$}
            \addlegendentry{$t_f=2.0\mathrm{E}{-3}$}
            	\label{plot_one} 
                
      \nextgroupplot[ legend cell align={left},
   						legend style={at={(axis cs:400,60)},anchor=south west, draw = none,line width=1.0pt},
                        xmin = 0, xmax = 1000,
                        xtick ={0,200,...,1000},
                        xlabel near ticks, ylabel near ticks,
                        ytick ={0,20,...,120},
                        ymin = 0, ymax = 120]

                        \addplot[red] table [x=x ,y=tf1] {plot/con_sigma_0.06.txt};
                        \addplot[c3] table [x=x ,y=tf2] {plot/con_sigma_0.06.txt};
                        \addplot[blue!70] table [x=x ,y=tf3] {plot/con_sigma_0.06.txt};
                        \addplot[black] table [x=x ,y=tf4] {plot/con_sigma_0.06.txt};
                        
            \addlegendimage{/pgfplots/refstyle=plot_two}
            \addlegendentry{$t_f=2.0\mathrm{E}{0}$}
            \addlegendentry{$t_f=2.0\mathrm{E}{-1}$}
            \addlegendentry{$t_f=2.0\mathrm{E}{-2}$}
            \addlegendentry{$t_f=2.0\mathrm{E}{-3}$}
            	\label{plot_two}

	\end{groupplot}
	
		\end{tikzpicture}
	\caption{Eigenvalues of the matrices for different values of fault transmissibility: (left) the mixed method is simulated in the mesh with $h=1.0\mathrm{E}{-1}$ and dof=$2210$, (right) the new method is simulated in the mesh with $h=5.0\mathrm{E}{-2}$, $\epsilon = 3h_f$, and dof=$3279$.}
	\label{eigen}
\end{figure}
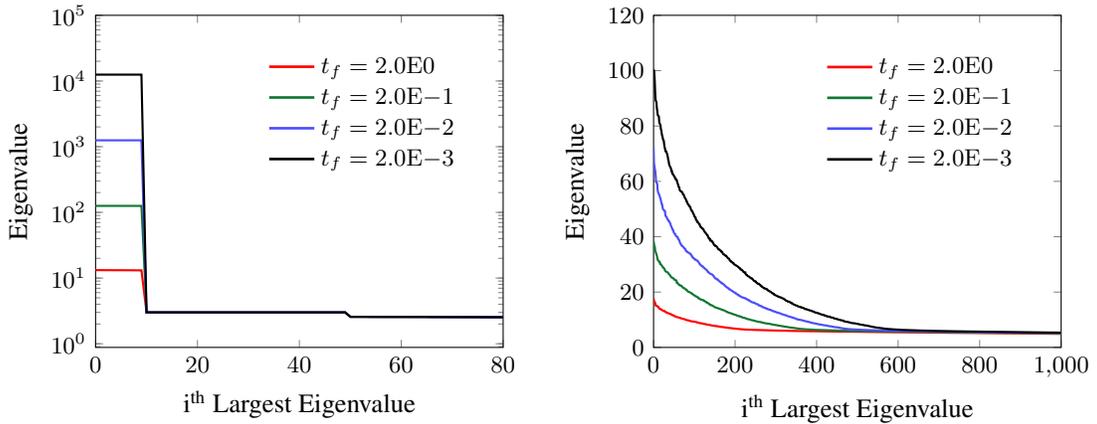

We also reported the effects of the transmissibility values on the matrices' eigenvalues of both the two methods in \autoref{eigen}. It can be seen that $t_f$ has large effect on the largest eigenvalue, and thus the condition number of the matrices in both methods. Moreover, the largest eigenvalue scales as $t_f$ in the mixed method, and scales as log$t_f$ in the new method.

\begin{figure}[t!]
\centering
\begin{tikzpicture}[every plot/.append style={thick}]
\pgfplotsset{
scale only axis,
every axis legend/.append style={
legend columns=4
}
}
\begin{groupplot}[
footnotesize,
width=6cm,
height=4cm,
group style={
group name=case1, group size=2 by 2,
horizontal sep=2cm, vertical sep=1.5cm,
},
xlabel near ticks, ylabel near ticks
]
\nextgroupplot[legend cell align={left},
legend style={nodes={scale=1.2, transform shape},
anchor=south west,
draw = none,line width=1.0pt},
legend to name={Legend1},
xmin = 0, xmax = 2,
xtick ={0,0.4,...,2},
xlabel=$x$,
ylabel=$p$,
ytick ={0,0.2,...,1},
ymin = 0, ymax = 1]

\addplot[red] table [x=x ,y=y] {plot/2D_px_05_exact.txt};
\addplot[c4, dotted] table [x=x ,y=y] {plot/2D_px_4_3.0_05.txt};
\addplot[blue!70, loosely dashed] table [x=x ,y=y] {plot/2D_px_4_4.0_05.txt};
\addplot[black, densely dashed] table [x=x ,y=y] {plot/2D_px_4_5.0_05.txt};

\addlegendimage{/pgfplots/refstyle=plot_one}
\addlegendentry{\scriptsize{Ground truth}}
\addlegendentry{\scriptsize{$\epsilon = 3h_f$}}
\addlegendentry{\scriptsize{$\epsilon = 4h_f$}}
\addlegendentry{\scriptsize{$\epsilon = 5h_f$}}
\label{plot_one}

\nextgroupplot[ legend cell align={left},
legend style={at={(axis cs:1.1,0.5)},anchor=south west, draw = none,line width=1.0pt},
xmin = 0, xmax = 2,
xtick ={0,0.4,...,2},
xlabel=$x$,
ylabel=$u_n$,
xlabel near ticks, ylabel near ticks,
ytick ={0,0.2,...,1},
ymin = 0, ymax = 1]

\addplot[red] table [x=x ,y=y] {plot/2D_ux_05_exact.txt};
\addplot[c4, dotted] table [x=x ,y=y] {plot/2D_ux_4_3.0_05.txt};
\addplot[blue!70, loosely dashed] table [x=x ,y=y] {plot/2D_ux_4_4.0_05.txt};
\addplot[black, densely dashed] table [x=x ,y=y] {plot/2D_ux_4_5.0_05.txt};

\nextgroupplot[legend cell align={left},
legend style={nodes={scale=1.2, transform shape},
anchor=south west,
draw = none,line width=1.0pt},
legend to name={Legend2},
xmin = 0, xmax = 2,
xtick ={0,0.4,...,2},
xlabel=$x$,
ylabel=$p$,
ytick ={0,0.2,...,1},
ymin = 0, ymax = 1]

\addplot[red] table [x=x ,y=y] {plot/2D_p_exact_05_002.txt};
\addplot[c4, dotted] table [x=x ,y=y] {plot/2D_px_45_1.0_05.txt};
\addplot[blue!70, loosely dashed] table [x=x ,y=y] {plot/2D_px_45_2.0_05.txt};
\addplot[black, densely dashed] table [x=x ,y=y] {plot/2D_px_45_3.0_05.txt};

\addlegendimage{/pgfplots/refstyle=plot_one}
\addlegendentry{\scriptsize{Ground truth}}
\addlegendentry{\scriptsize{$\epsilon = 1h_f$}}
\addlegendentry{\scriptsize{$\epsilon = 2h_f$}}
\addlegendentry{\scriptsize{$\epsilon = 3h_f$}}
\label{plot_one}

\nextgroupplot[ legend cell align={left},
legend style={at={(axis cs:1.1,0.5)},anchor=south west, draw = none,line width=1.0pt},
xmin = 0, xmax = 2,
xtick ={0,0.4,...,2},
xlabel=$x$,
ylabel=$u_n$,
xlabel near ticks, ylabel near ticks,
ytick ={0,0.2,...,1},
ymin = 0, ymax = 1]

\addplot[red] table [x=x ,y=y] {plot/2D_ux_exact_05_002.txt};
\addplot[c4, dotted] table [x=x ,y=y] {plot/2D_ux_45_1.0_05.txt};
\addplot[blue!70, loosely dashed] table [x=x ,y=y] {plot/2D_ux_45_2.0_05.txt};
\addplot[black, densely dashed] table [x=x ,y=y] {plot/2D_ux_45_3.0_05.txt};

\end{groupplot}
\path (2d c1r1.north) -- node[above right = 0.5cm and -2cm,font=\fontsize{12pt}{5pt}\selectfont]{\ref{Legend1}} (2d c2r1.north);
\path (2d c1r2.north) -- node[above right =-0.3cm and -2cm,font=\fontsize{12pt}{5pt}\selectfont]{\ref{Legend2}} (2d c2r2.north);

\end{tikzpicture}
\caption{Pressure (left) and normal component of velocity (right) solutions from the ground truth and the new method of mesh with $h =5.0\mathrm{E}{-3}$ and $\epsilon = 3h_f$ along the line $\{(x,0.5): 0 \leq x \leq L_x\}$ in the 2D test when $t_f = 2.0$ (upper) and $t_f = 0.02$ (lower).}
\label{2d_plot}
\end{figure}
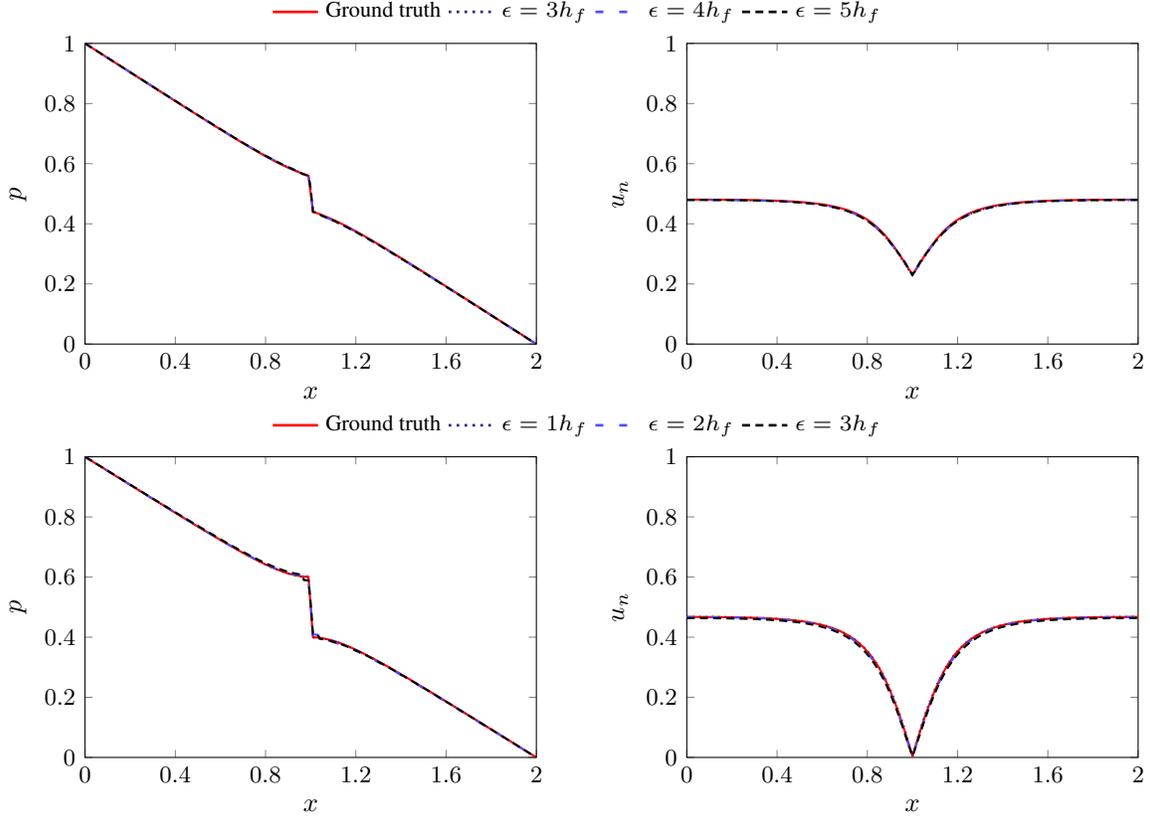

\clearpage
\subsection{3D test}\label{3D}
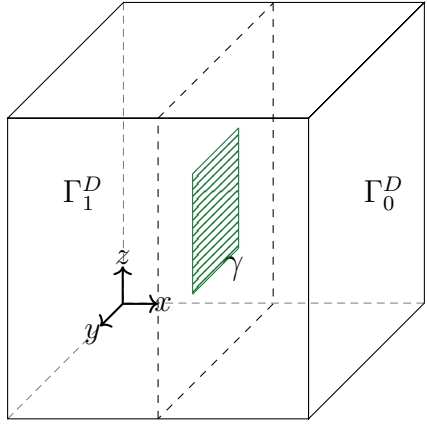
\begin{figure}[t!]
    \centering
    \begin{tikzpicture}
    \pic at (1,0) {annotated cuboid={width=40, height=40, depth=40}};
    \draw[dashed] (-1.0,0.) -- (0.53,1.53);
    \draw[dashed] (0.53,1.53) -- (0.53,-2.47);
    \draw[dashed] (0.53,-2.47) -- (-1, -4);
    \draw[dashed] (-1, -4) -- (-1,0);
    \node[] at (0, -2) {\Large$\gamma$};
    \node[] at (-2, -1) {\large$\Gamma_1^D$};
    \node[] at (2, -1) {\large$\Gamma_0^D$};

    \draw[thick, ->] (-1.47, -2.47) to (-1,-2.47);
    \draw[thick, ->] (-1.47, -2.47) to (-1.47,-1.97);
    \draw[thick, ->] (-1.47, -2.47) to (-1.77,-2.77);
    \node[] at (-0.93,-2.47)   (a) {\large$x$};
    \node[] at (-1.87,-2.87)   (a) {\large$y$};
    \node[] at (-1.47,-1.85)   (a) {\large$z$};
    
    \draw[c3] (-0.541, -0.741) -- (0.071, -0.129);
    \draw[c3] (0.071, -0.129) -- (0.071, -1.729);
    \draw[c3] (0.071, -1.729) -- (-0.541, -2.341);
    \draw[c3] (-0.541, -2.341) -- (-0.541, -0.741);
    \fill[pattern=flexible north east lines, pattern density = 3pt, pattern line width=.5pt,pattern color=c3] (-0.541, -0.741) -- (0.071, -0.129) -- (0.071, -1.729) -- (-0.541, -2.341) --  (-0.541, -0.741) ;
 
    \end{tikzpicture}
    \caption{Configuration of the 3D test. }
    \label{3d_setup}
\end{figure}

\paragraph{Setup.} The 3D test case is described in \autoref{3d_setup} where the domain $\Omega$ is a unit cube. A square-shaped fault $\gamma$ is located at the center of the plane $\{(0.5,x,y):0\leq y,z\leq 1\}$ in the domain. The length of each edge on $\gamma$ is 0.4 and a small transmissibility $t_f = $0.0005 is considered. We define different part of boundaries by $\Gamma_1^D = \{(0,y, z): 0 \leq y,z \leq 1 \}$, $\Gamma_0^D = \{(1,y,z): 0 \leq y,z \leq 1\}$, and $\Gamma^N = \partial \Omega \backslash (\Gamma_1^D \cup \Gamma_0^D) $. The boundary conditions are defined by
\begin{align*}
\begin{cases}
    p = 1 \qquad \text{on $\Gamma_1^D$},\\
    p = 0  \qquad \text{on $\Gamma_0^D$},\\
    \nabla \mathbf{u}\cdot \mathbf{\nu} = 0 \qquad \text{on $\Gamma^N$}.
\end{cases}
\end{align*}

The meshes are defined in the following way that is similar to the 2D test. We denote the mesh size on the domain boundary $\partial \Omega$ by $h$, the mesh size on the boundary of the subdomain $\partial \Omega_s$ by $h_s$, and mesh size on the fault $\gamma$ by $h_f$. The three mesh sizes have the relation of $h:h_s:h_f = 5:2:2$ for different case. Moreover, we set the distance between subdomain boundary to the fault as $10h_f$. 

The mixed method and the new method are both conducted at five different meshes for comparison. Furthermore, three values of $\epsilon$ are considered for the new method: $1h_f$, $2h_f$ and $3h_f$.

\begin{figure}[t!]
\centering
\begin{tikzpicture}[every plot/.append style={thick}]
\pgfplotsset{
scale only axis,
every axis legend/.append style={
legend columns=4
}
}
\begin{groupplot}[
xlabel near ticks, ylabel near ticks,
width=4cm,
height=4cm,
group style={
group name=2d, group size=3 by 2,
horizontal sep=.5cm, vertical sep=1.2cm,
y descriptions at=edge left,
}
]
\nextgroupplot[legend cell align={left},
legend style={nodes={scale=1.2, transform shape},
anchor=south west,
draw = none,line width=1.0pt},
legend to name={CommonLegend},
xlabel=$h$,
xmin = 0.005, xmax = 0.2,
xmode = log,
ymax = 3.e-2,ymin = 1.e-3,
ylabel=$\Vert e_p^h \Vert_{L^2}$,
ymode=log,
log basis y={10},
xlabel near ticks, ylabel near ticks,
grid = both]
\addplot[red,mark= square, mark options={solid}] table [x=h ,y=ep]
{fig/3D_mix_error.txt};
\addplot[c3,mark= o,mark options={solid}] table [x=h ,y=ep_1]
{fig/3D_con_error.txt};
\addplot[blue!70,mark= asterisk, mark options={solid}] table [x=h ,y=ep_2]
{fig/3D_con_error.txt};
\addplot[black,mark= diamond,mark options={solid}] table [x=h ,y=ep_3]
{fig/3D_con_error.txt};
\draw[black,thick] (axis cs:0.01,0.005) -- node[left= 0.05cm] {$1$} (axis cs:0.01,0.01);
\draw[black,thick] (axis cs:0.01,0.01) -- node[above= 0.05cm] {$2$} (axis cs:0.04,0.01);
\draw[black,thick] (axis cs:0.04,0.01) --  (axis cs:0.01,0.005);

\draw[black,thick] (axis cs:0.05,0.002) -- node[below= 0.05cm] {$1$} (axis cs:0.1,0.002);
\draw[black,thick] (axis cs:0.1,0.002) -- node[right= 0.05cm] {$1$} (axis cs:0.1,0.004);
\draw[black,thick] (axis cs:0.1,0.004) --  (axis cs:0.05,0.002);

\addlegendimage{/pgfplots/refstyle=plot_one}
\addlegendentry{\scriptsize{mixed method}}
\addlegendentry{\scriptsize{$\epsilon = 1h_f$}}
\addlegendentry{\scriptsize{$\epsilon = 2h_f$}}
\addlegendentry{\scriptsize{$\epsilon = 3h_f$}}

\nextgroupplot[ legend cell align={left},
xlabel near ticks, ylabel near ticks,
legend style={at={(axis cs:10000,0.03)},nodes={scale=0.8, transform shape},anchor=south west, draw = none,line width=1.0pt},
xlabel=dof,
xmin = 10000, xmax = 10000000,
xmode = log,
ymax = 3.e-2,ymin = 1.e-3,
ymode=log,
log basis y={10},
grid = both]
\addplot[red,mark= square, mark options={solid}] table [x=dof ,y=ep]
{fig/3D_mix_error.txt};
\addplot[c3,mark= o,mark options={solid}] table [x=dof,y=ep_1]
{fig/3D_con_error.txt};
\addplot[blue!70,mark= asterisk, mark options={solid}] table [x=dof ,y=ep_2]
{fig/3D_con_error.txt};
\addplot[black,mark= diamond,mark options={solid}] table [x=dof ,y=ep_3]
{fig/3D_con_error.txt};

\nextgroupplot[ legend cell align={left},
xlabel near ticks, ylabel near ticks,
legend style={at={(axis cs:0.1,0.03)},nodes={scale=0.8, transform shape},anchor=south west, draw = none,line width=1.0pt},
xlabel=time,
xmin = 0.1, xmax = 1000,
xmode = log,
ymax = 3.e-2,ymin = 1.e-3,
ymode=log,
log basis y={10},
grid = both]
\addplot[red,mark= square, mark options={solid}] table [x=t ,y=ep]
{fig/3D_mix_error.txt};
\addplot[c3,mark= o,mark options={solid}] table [x=t_1,y=ep_1]
{fig/3D_con_error.txt};
\addplot[blue!70,mark= asterisk, mark options={solid}] table [x=t_2 ,y=ep_2]
{fig/3D_con_error.txt};
\addplot[black,mark= diamond,mark options={solid}] table [x=t_3 ,y=ep_3]
{fig/3D_con_error.txt};

\nextgroupplot[ legend cell align={left},
legend style={at={(axis cs:0.005,0.08)},nodes={scale=0.8, transform shape},anchor=south west, draw = none,line width=1.0pt},
xlabel=$h$,
xmin = 0.005, xmax = 0.2,
xmode = log,
ymax = 1.e-1,ymin = 1.e-2,
ylabel=$\Vert \mathbf{e_u^h} \Vert_{L^2}$,
ymode=log,
log basis y={10},
xlabel near ticks, ylabel near ticks,
grid = both]
\addplot[red,mark= square, mark options={solid}] table [x=h ,y=eu]
{fig/3D_mix_error.txt};
\addplot[c3,mark= o,mark options={solid}] table [x=h ,y=eu_1]
{fig/3D_con_error.txt};
\addplot[blue!70,mark= asterisk, mark options={solid}] table [x=h ,y=eu_2]
{fig/3D_con_error.txt};
\addplot[black,mark= diamond,mark options={solid}] table [x=h ,y=eu_3]
{fig/3D_con_error.txt};
\draw[black,thick] (axis cs:0.01,0.06) -- node[above= 0.05cm] {$2$} (axis cs:0.04, 0.06);
\draw[black,thick] (axis cs:0.01,0.06) -- node[left= 0.05cm] {$1$} (axis cs:0.01,0.03);
\draw[black,thick] (axis cs:0.01,0.03) --  (axis cs:0.04,0.06);

\nextgroupplot[ legend cell align={left},
xlabel near ticks, ylabel near ticks,
legend style={at={(axis cs:10000,0.08)},nodes={scale=0.8, transform shape},anchor=south west, draw = none,line width=1.0pt},
xlabel=dof,
xmin = 10000, xmax = 10000000,
xmode = log,
ymax = 2.e-1,ymin = 5.e-3,
ymode=log,
log basis y={10},
grid = both]
\addplot[red,mark= square, mark options={solid}] table [x=dof ,y=eu]
{fig/3D_mix_error.txt};
\addplot[c3,mark= o,mark options={solid}] table [x=dof,y=eu_1]
{fig/3D_con_error.txt};
\addplot[blue!70,mark= asterisk, mark options={solid}] table [x=dof ,y=eu_2]
{fig/3D_con_error.txt};
\addplot[black,mark= diamond,mark options={solid}] table [x=dof ,y=eu_3]
{fig/3D_con_error.txt};

\nextgroupplot[ legend cell align={left},
xlabel near ticks, ylabel near ticks,
legend style={at={(axis cs:0.1,0.08)},nodes={scale=0.8, transform shape},anchor=south west, draw = none,line width=1.0pt},
xlabel=time,
xmin = 0.1, xmax = 1000,
xmode = log,
ymax = 2.e-1,ymin = 5.e-3,
ymode=log,
log basis y={10},
grid = both]
\addplot[red,mark= square, mark options={solid}] table [x=t ,y=eu]
{fig/3D_mix_error.txt};
\addplot[c3,mark= o,mark options={solid}] table [x=t_1,y=eu_1]
{fig/3D_con_error.txt};
\addplot[blue!70,mark= asterisk, mark options={solid}] table [x=t_2 ,y=eu_2]
{fig/3D_con_error.txt};
\addplot[black,mark= diamond,mark options={solid}] table [x=t_3 ,y=eu_3]
{fig/3D_con_error.txt};

\end{groupplot}
\path (2d c2r1.north west) -- node[above=0.2cm,font=\fontsize{12pt}{5pt}\selectfont]{\ref{CommonLegend}} (2d c2r1.north east);
\end{tikzpicture}
\caption{$L^2$ errors of pressure and velocity for the new method and the mixed method in the 3D test.}
\label{3d_convergence}
\end{figure}
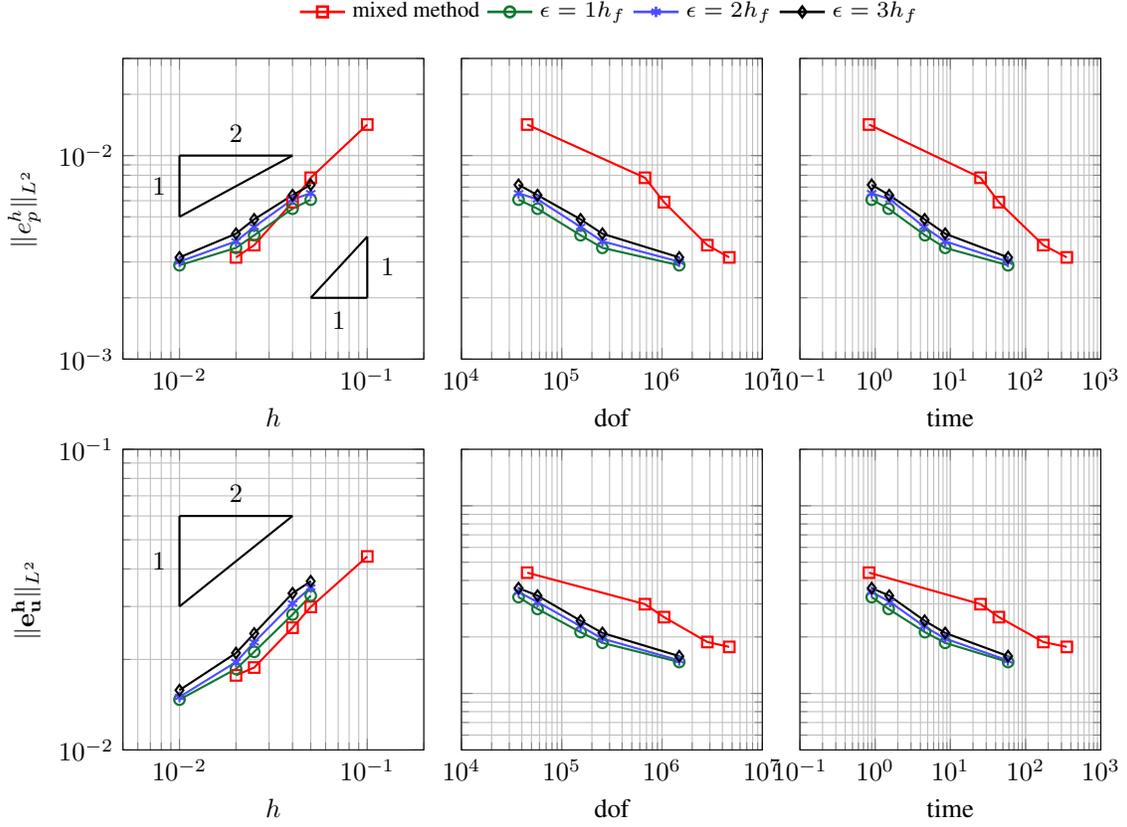

\paragraph{Results.} The simulation results of the mixed method and the new method are summarized in \autoref{3d_table_mix} and \autoref{3d_table_con}, respectively. The $L^2$ pressure error and velocity error are plotted as functions of $h$, dof and simulation time for different cases in \autoref{3d_convergence}. From these plots, following conclusions could be derived: first, the pressure solution converges as $\mathcal{O}(h^{\frac{1}{2}})$ in the new method, and converges as $\mathcal{O}(h)$ in the mixed method. The pressure $L^2$ error of the new method could be lower then the error of the mixed method at some coarse mesh size. The velocity solution converges as $\mathcal{O}(h^{\frac{1}{2}})$ toward the ground truth solution for both of the methods. Second, at given dof or simulation time, the new method generates lower $L^2$ error than the mixed method for both the pressure and velocity solutions in the tests shown here. Third, we can observe that test with lower value $\epsilon$ has lower $L^2$ error in the new method. Among the tests conducted with different $\epsilon$'s here, $\epsilon = 1h_f$ generates the optimal solution.

\begin{figure}[t!]
\centering
    \begin{tikzpicture}
    	\node[inner sep=0pt] (slice1) at (2, 0)
    		{\includegraphics[width=5.5cm,height = 4.5cm]{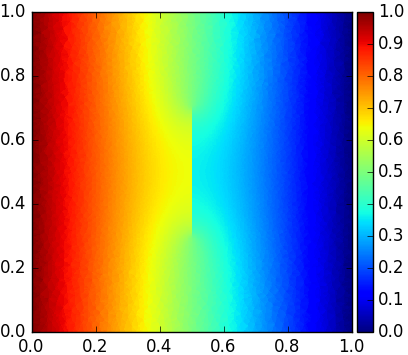}};
    	  \node[inner sep=0pt] (slice2) at (9.5, 0)
    		{\includegraphics[width=5.5cm,height = 4.5cm]{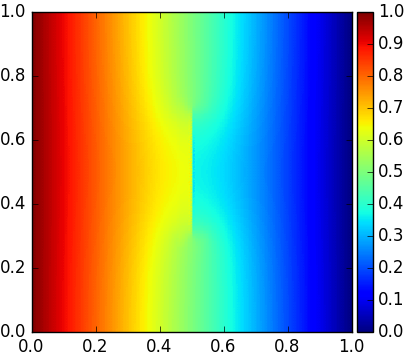}};
    	  \node[inner sep=0pt] (slice3) at (2, -4.7)
    		{\includegraphics[width=5.5cm,height = 4.5cm]{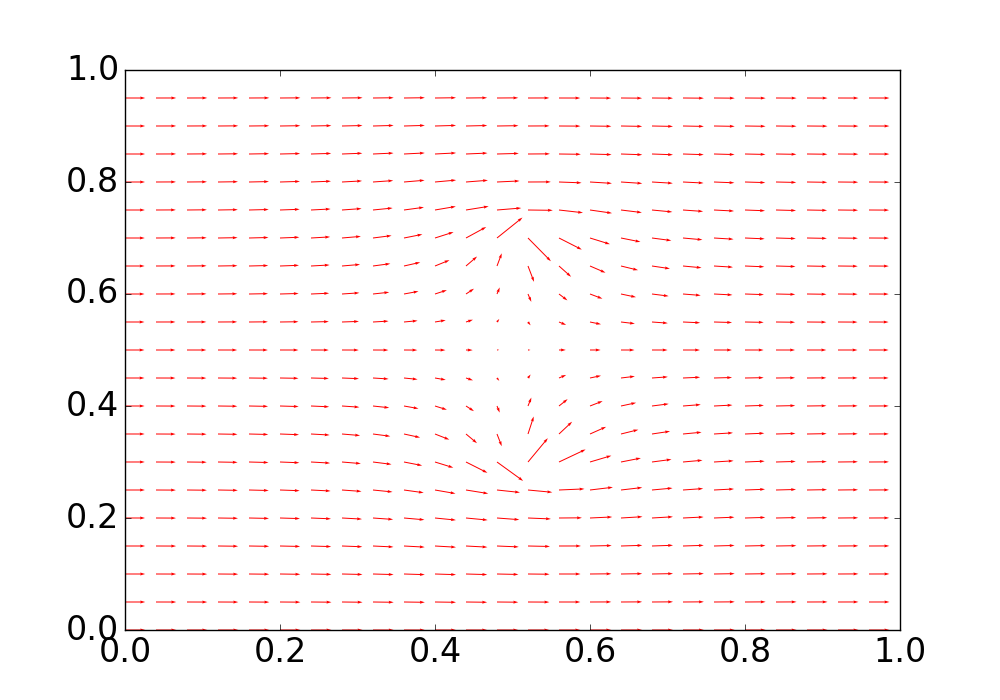}};
    	\node[inner sep=0pt] (slice3) at (9.5, -4.7)
    		{\includegraphics[width=5.5cm,height = 4.5cm]{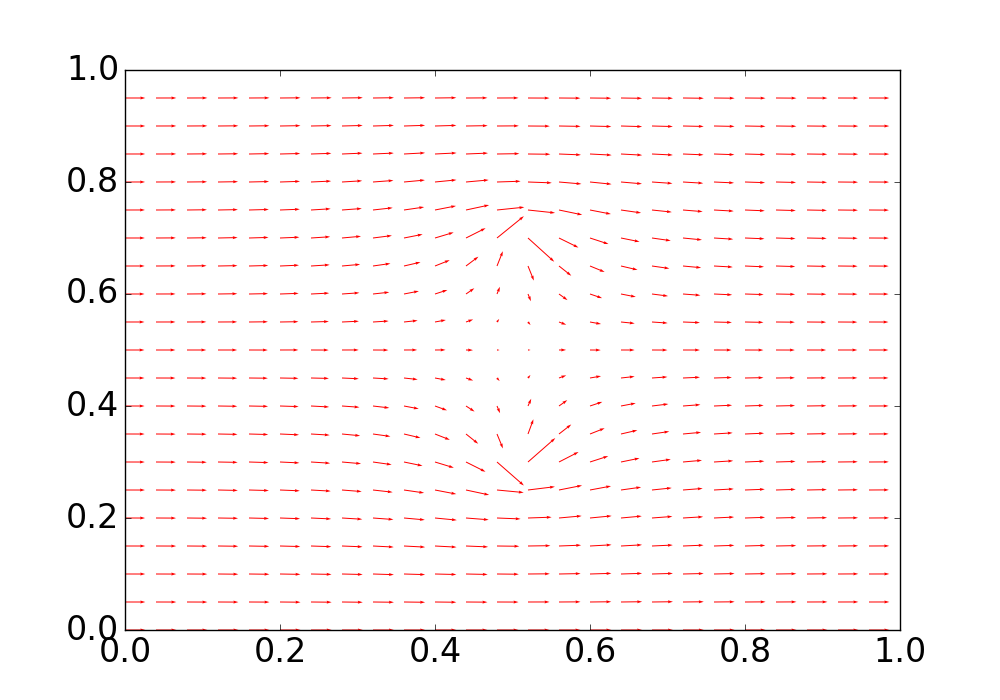}};
    \end{tikzpicture}
    \caption{Pressure and velocity solutions on the plane $\{(x,y,0.5):0\leq x,y \leq 1 \}$: (left) ground truth solutions, (right) solutions of the new method with $h =1.0\mathrm{E}{-2}$ and $\epsilon = 1h_f$ in the 3D test.}
    \label{3d_p}
\end{figure}

\begin{figure}[t!]
\centering
\begin{tikzpicture}[every plot/.append style={thick}]
\pgfplotsset{
scale only axis,
every axis legend/.append style={
legend columns=4
}
}
\begin{groupplot}[
footnotesize,
width=5cm,
height=3.7cm,
group style={
group name=case1, group size=2 by 2,
horizontal sep=2cm, vertical sep=1.2cm,
},
xlabel near ticks, ylabel near ticks
]
\nextgroupplot[legend cell align={left},
legend style={nodes={scale=1.2, transform shape},
anchor=south west,
draw = none,line width=1.0pt},
legend to name={Legend1},
xmin = 0, xmax = 1,
xtick ={0,0.2,...,1},
xlabel=$x$,
ylabel=$p$,
ytick ={0,0.2,...,1},
ymin = 0, ymax = 1]

\addplot[red] table [x=x ,y=y] {plot/mix_p_0.4.txt};
\addplot[c4, dotted] table [x=x ,y=y] {plot/con_p_1_0.4.txt};
\addplot[blue!70, loosely dashed] table [x=x ,y=y] {plot/con_p_2_0.4.txt};
\addplot[black, densely dashed] table [x=x ,y=y] {plot/con_p_3_0.4.txt};

\addlegendimage{/pgfplots/refstyle=plot_one}
\addlegendentry{\scriptsize{Ground truth}}
\addlegendentry{\scriptsize{$\epsilon = 1h_f$}}
\addlegendentry{\scriptsize{$\epsilon = 2h_f$}}
\addlegendentry{\scriptsize{$\epsilon = 3h_f$}}

\nextgroupplot[ legend cell align={left},
legend style={at={(axis cs:0.55,.95)},anchor=south west, draw = none,line width=1.0pt},
xmin = 0, xmax = 1,
xtick ={0,0.2,...,1},
xlabel=$x$,
ylabel=$u_n$,
xlabel near ticks, ylabel near ticks,
ytick ={0,0.4,...,2},
ymin = 0, ymax = 2]

\addplot[red] table [x=x ,y=y] {plot/mix_u_0.4.txt};
\addplot[c4, dotted] table [x=x ,y=y] {plot/con_u_1_0.4.txt};
\addplot[blue!70, loosely dashed] table [x=x ,y=y] {plot/con_u_2_0.4.txt};
\addplot[black, densely dashed] table [x=x ,y=y] {plot/con_u_3_0.4.txt};

\nextgroupplot[ legend cell align={left},
legend style={at={(axis cs:0.55,0.5)},anchor=south west, draw = none,line width=1.0pt},
xmin = 0, xmax = 1,
xtick ={0,0.2,...,1},
xlabel=$x$,
ylabel=$p$,
ytick ={0,0.2,...,1},
ymin = 0, ymax = 1]

\addplot[red] table [x=x ,y=y] {plot/mix_p_0.5.txt};
\addplot[c4, dotted] table [x=x ,y=y] {plot/con_p_1_0.5.txt};
\addplot[blue!70, loosely dashed] table [x=x ,y=y] {plot/con_p_2_0.5.txt};
\addplot[black, densely dashed] table [x=x ,y=y] {plot/con_p_3_0.5.txt};
\nextgroupplot[ legend cell align={left},
legend style={at={(axis cs:0.55,0.95)},anchor=south west, draw = none,line width=1.0pt},
xmin = 0, xmax = 1,
xtick ={0,0.2,...,1},
xlabel=$x$,
ylabel=$u_n$,
xlabel near ticks, ylabel near ticks,
ytick ={0,0.4,...,2},
ymin = 0, ymax = 2]

\addplot[red] table [x=x ,y=y] {plot/mix_u_0.5.txt};
\addplot[c4, dotted] table [x=x ,y=y] {plot/con_u_1_0.5.txt};
\addplot[blue!70, loosely dashed] table [x=x ,y=y] {plot/con_u_2_0.5.txt};
\addplot[black, densely dashed] table [x=x ,y=y] {plot/con_u_3_0.5.txt};

\end{groupplot}
\path (2d c1r1.north) -- node[above right = -0.1cm and -2cm,font=\fontsize{12pt}{5pt}\selectfont]{\ref{Legend1}} (2d c2r1.north);
\end{tikzpicture}
\caption{Pressure and normal component of velocity solutions along the line $\{(x,0.4,0.5): 0 \leq x \leq 1\}$ (upper) and along the line $\{(x,0.5,0.5): 0 \leq x \leq 1\}$ (lower): the new method is conducted at $h=1.0\mathrm{E}{-2}$ in the 3D test.}
\label{3d_plot}
\end{figure}
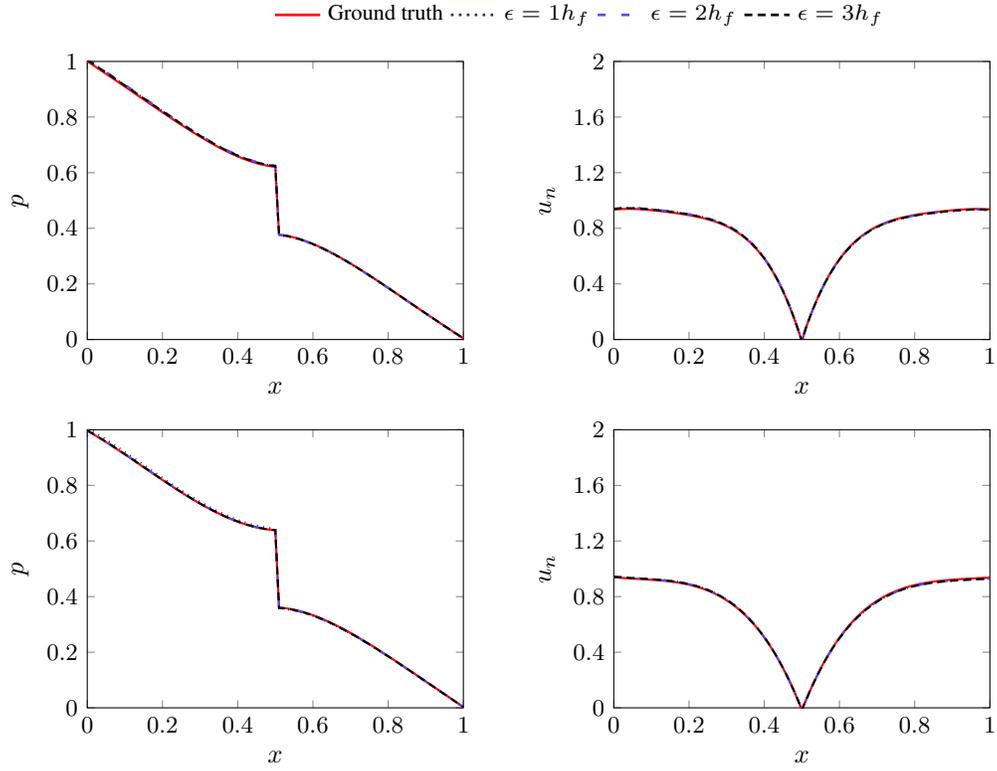

The pressure and velocity solutions of the ground truth test on the plane $\{(x,y,0.5):0\leq x,y\leq1 \}$ is shown in \autoref{3d_p}. For comparison, we also plot the solutions of the new method with mesh size $h=1.0\mathrm{E}{-2}$ and $\epsilon = 1h_f$. Moreover, the pressure and normal component of velocity solutions along the line $\{(x,0.4,0.5): 0 \leq x \leq 1 \}$ and along the line $\{(x,0.5,0.5): 0 \leq x \leq 1 \}$ are plotted in \autoref{3d_plot}. It can be seen that the solutions generated by the new method match the ground truth solutions very well. We can also observe that the low-transmissibility fault in this test behaves like a barrier and generates high pressure jump and near-zero normal component of velocity along the fault. 

\setlength{\extrarowheight}{5pt}
\begin{table}[t]
	\scriptsize
	\centering
    \caption{Summary of simulation results for the mixed method in the 3D test. Number in the bracket of "time" column represents the preconditioning time of ILU1 preconditioner. "*" represents the ground truth test. }
	\label{3d_table_mix}
    \begin{tabular*}{\textwidth}{@{\extracolsep{7pt}}
                                                         P{0.18\textwidth-2\tabcolsep}
                                                         P{0.18\textwidth-2\tabcolsep}
                                 P{0.18\textwidth-2\tabcolsep}
                                 P{0.18\textwidth-2\tabcolsep}
                                 P{0.18\textwidth-2\tabcolsep}}
    \toprule %
    h  &dof  &time &$\Vert e_p^h \Vert_{L^2}$ &$\Vert \mathbf{e_u^h} \Vert_{L^2}$  \\
    \toprule
1.0E-01   &45273   &0.83 (0.2)   &1.42E-02   &4.40E-02   \\
5.0E-02   &676033   &25.19 (5.5)  &7.78E-03   &2.99E-02   \\
4.0E-02   &1047652   &44.16 (8.6)  &5.90E-03   &2.55E-02   \\
2.5E-02   &2829944   &171.94 (20.2)  &3.63E-03   &1.88E-02   \\
2.0E-02   &4684445   &353.27 (38.4)  &3.16E-03   &1.77E-02   \\
1.0E-02\textsuperscript{*}   &14447330   &3186.33   &-   &- \\
    \bottomrule
 	\end{tabular*}
\end{table}

\begin{table}[t]
	\scriptsize
	\centering
    \caption{Summary of simulation results of the new method in the 3D test.The first number in "time" column represents the CPU time used for solving the approximate pressure and velocity of the whole domain, the second number is the time for the subdomain problem. Number in the bracket of "time" column represents the preconditioning time of ILU1 preconditioner.}
	\label{3d_table_con}
    \begin{tabular*}{\textwidth}{@{\extracolsep{2pt}}
    							 P{0.085\textwidth-2\tabcolsep}
    							 P{0.08\textwidth-2\tabcolsep}
                                 P{0.12\textwidth-2\tabcolsep}
                                 P{0.085\textwidth-2\tabcolsep}
                                 P{0.085\textwidth-2\tabcolsep}
                                 P{0.086\textwidth-2\tabcolsep}
                                 P{0.085\textwidth-2\tabcolsep}
                                 P{0.085\textwidth-2\tabcolsep}
                                 P{0.086\textwidth-2\tabcolsep}
                                 P{0.085\textwidth-2\tabcolsep}
                                 P{0.085\textwidth-2\tabcolsep}}
    \toprule %
    \multirow{2}{*}{h} &\multirow{2}{*}{dof} &\multicolumn{3}{c}{$\epsilon = 1h_f$} &\multicolumn{3}{c}{$\epsilon = 2h_f$} &\multicolumn{3}{c}{$\epsilon = 3h_f$}\\
    \cline{3-5} \cline{6-8} \cline{9-11}
     & &time &$\Vert e_p^h\Vert_{L^2}$ &$\Vert \mathbf{e_{\mathbf{u}}^h} \Vert_{L^2}$ &time &$\Vert e_p^h\Vert_{L^2}$ &$\Vert \mathbf{e_{\mathbf{u}}^h} \Vert_{L^2}$ &time &$\Vert e_p^h\Vert_{L^2}$ &$\Vert \mathbf{e_{\mathbf{u}}^h} \Vert_{L^2}$\\
    \toprule
    5.0E-02   &37115   &0.4+0.5 (0.1)   &6.1E-03   &3.3E-02   &0.4+0.5   &6.5E-03   &3.5E-02   &0.4+0.5   &7.2E-03   &3.6E-02  \\ 
4.0E-02   &57410   &0.6+0.9 (0.2)   &5.5E-03   &2.8E-02   &0.7+0.9  &6.1E-03   &3.1E-02   &0.7+0.9   &6.4E-03   &3.3E-02   \\
2.5E-02   &154093   &2.8+1.8 (0.5)   &4.1E-03   &2.1E-02   &2.7+1.9   &4.5E-03   &2.3E-02   &2.8+1.8   &4.9E-03   &2.4E-02   \\
2.0E-02   &254306   &5.5+3.1 (0.9)  &3.5E-03   &1.9E-02   &5.5+3.1   &3.8E-03   &2.0E-02   &5.5+3.1   &4.1E-03   &2.1E-02   \\
1.0E-02   &1485031   &37.0+21.6 (5.1)   &2.9E-03   &1.5E-02   &37.0+21.3   &3.0E-03   &1.5E-02   &37.0+21.5   &3.2E-03   &1.6E-02\\  
    \bottomrule
 	\end{tabular*}
\end{table}

\section{Conclusion} \label{conclusion}
We presented in this paper a new method to solve pressure and velocity of flow problem in porous media with low-permeability faults. We defined pressures as approximations to the pressure in the original boundary value problem. Then we decoupled the approximate pressure from velocity and obtained new formulations for the pressure that could be solved by continuous Galerkin finite element method. To correct the approximations around the fault and obtain the final solutions, we additionally solved small problem associated with subdomain around the fault using the mixed finite element method.
We conducted numerical tests in both two-dimensional domain and three-dimensional domains to investigate the accuracy of the new method. Error results showed that the convergence rate for pressure solution of the new method depends on the fault transmissibility, but velocity solution generally has convergence rate of $\mathcal{O}(h^\frac{1}{2})$. Compared to the mixed method, the new method proposed in this paper can be faster at given pressure and velocity error in three-dimensional test.

\clearpage
\bibliographystyle{unsrt}  
\bibliography{references}

\clearpage
\begin{appendices}
\section{Formulation of \texorpdfstring{$-\Delta p$}{Lg}}\label{appendices}

\begin{figure}[t!]
\begin{center}
  \begin{tikzpicture}
        \path[font={\tiny}]
        (11 , -.8)   coordinate (A)
        (5.7   , -.8)   coordinate (B)
        (5.7   , 3)   coordinate (C)
        (11   , 3)   coordinate (D)
    ;
    \draw plot [smooth cycle, tension=0.8] coordinates {(A) (B) (C) (D)};
    \draw[black, line width = 0.2mm] (8.5,0) to (8.5,2.5);
    \draw[black, dashed, line width = 0.2mm] (8.5,-1.35) to (8.5, 0);
    \draw[black, dashed, line width = 0.2mm] (8.5, 2.5) to (8.5, 3.6);
        \node[] at (6.3,0.5)   (a) {$\Omega_1$};
        \node[] at (10.5,0.5)   (a) {$\Omega_2$};
        \node[] at (4.6,1.5)   (a) {$\Gamma_1$};
        \node[] at (12.3,1.5)   (a) {$\Gamma_2$};
        \node[] at (8.3, 1.5) (c) {$\gamma^+$};
        \node[] at (8.8,1.5) (d) {$\gamma^-$};
	\node[] at (8.8, 2.9) (e) {$\gamma^-_e$};
        \node[] at (8.3, 2.9) (e) {$\gamma^+_e$};
        	\node[] at (8.8, -0.5) (e) {$\gamma^-_e$};
        \node[] at (8.3, -0.5) (e) {$\gamma^+_e$};
        \node[] at (9.3,0.5) (d) {$\mathbf{n}$};
    \draw[thick, ->] (8.5, 0.5) to (9,0.5);
  \end{tikzpicture}
  \end{center}
  \caption{The domain $\Omega$ with an immersed fault $\gamma$: the domain is split into subdomains $\Omega_1$ and $\Omega_2$ along $\gamma$.}
  \label{domain_divide}
\end{figure}
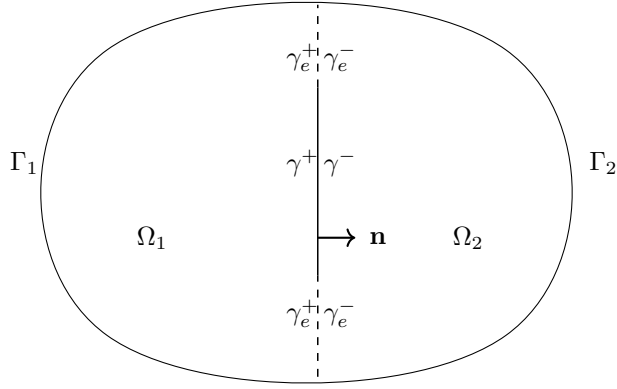

In this appendix we derive the differential equation for $-\Delta p$ of boundary value problem \eqref{BVP} by expressing the effects of faults as source term. For simplicity of discussion, we assume that $\perm$ is an identity matrix. Combine the first and second equations of \eqref{BVP_domain}, we can get
\begin{align}\label{deltap}
    -\Delta p = f \qquad \text{in $\Omega \backslash\gamma$},
\end{align}
and thus $-\Delta p \in L^2(\Subdomain)$. But there exits jumps for $p$ along $\gamma$, $-\Delta p$ in general is expected to be a distribution in $\mathcal{D}^\prime (\Omega)$. Here we are interested in the relationship between $-\Delta p$ and the jump condition.

We split the domain $\Omega$ into two subdomains $\Omega_1$ and $\Omega_2$ along $\gamma$ such that $\partial \Omega_1$, $\partial \Omega_2$ are Lipschitz (see \autoref{domain_divide}). Denote the extension part of $\gamma$ as $\gamma_e$ such that $\gamma \cup \gamma_e = \partial \Omega_1 \cap \partial \Omega_2$, and define $\Gamma_1 = \partial \Omega_1 \cap \partial \Omega$ and $\Gamma_2 = \partial \Omega_2 \cap \partial \Omega$, then we have
\begin{align*}
    \partial \Omega_1 = \Gamma_1 \cup \gamma^+ \cup \gamma^+_e, \qquad
    \partial \Omega_2 = \Gamma_2 \cup \gamma^- \cup \gamma^-_e.
\end{align*}

Take any test function $\phi$ from $\mathcal{D}(\Omega)$ and apply Green's second identity to $-\Delta p$, one can get
    \begin{align*}
	    -\langle \Delta p, \phi \rangle_{\Omega\backslash\gamma} = -\langle p ,\Delta \phi \rangle_{\Omega_1 \cup \Omega_2} \nolinebreak
	                                            + \langle p,  \nabla \phi \cdot \mathbf{\nu} \rangle_{\partial\Omega_1 \cup \partial\Omega_2} \nolinebreak
	                                            - \langle \nabla p \cdot \mathbf{\nu},  \phi \rangle_{\partial\Omega_1 \cup \partial\Omega_2}
    \end{align*}
where    
    \begin{align*}
        \langle p,  \nabla \phi \cdot \mathbf{\nu} \rangle_{\partial\Omega_1 \cup \partial\Omega_2} \nolinebreak
        &= \langle p,  \nabla \phi \cdot \mathbf{\nu} \rangle_{\partial \Omega \cup \gamma^+ \cup \gamma^- \cup \gamma_e^+ \cup \gamma_e^-} \\
        &= \langle p_D,  \nabla \phi \cdot \mathbf{\nu} \rangle_{\partial \Omega} \nolinebreak
        + \langle p|_{\gamma^+},  \nabla \phi \cdot \mathbf{n} \rangle_{\partial \gamma^+} \nolinebreak
        + \langle p|_{\gamma^-},  \nabla \phi \cdot \mathbf{-n} \rangle_{\partial \gamma^-} \\
        &= \langle p_D,  \nabla \phi \cdot \mathbf{\nu} \rangle_{\partial \Omega} \nolinebreak
         + \langle \pjump,  \nabla \phi \cdot \mathbf{n} \rangle_{\partial \gamma}
    \end{align*}
 and   
    \begin{align*}
        \langle \nabla p \cdot \mathbf{\nu},  \phi \rangle_{\partial\Omega_1 \cup \partial\Omega_2} \nolinebreak
         &= \langle \nabla p \cdot \mathbf{\nu},  \phi \rangle_{\partial \Omega }\nolinebreak
         + \langle \nabla p \cdot \mathbf{\nu},  \phi \rangle_{\gamma^+ \cup \gamma^- \cup \gamma_e^+ \cup \gamma_e^-} \\
         &= 0
    \end{align*}
since $\phi$ has compact support in $\Omega$, $\nabla p \cdot \mathbf{\nu}$ and $\phi$ are continuous along $\gamma$ and $\gamma_e$. Together with \eqref{deltap}, we can get the following equation:
    \begin{align}\label{L_p_H2}
        \begin{split}
        \langle f, \phi \rangle_{\Omega} &= -\langle \Delta p, \phi \rangle_{\Omega\backslash\gamma}\\
            &= -\langle p ,\Delta \phi \rangle_{\Omega} \nolinebreak
             + \langle p_D,  \nabla \phi \cdot \mathbf{\nu} \rangle_{\partial \Omega} \nolinebreak
            + \langle \pjump,  \nabla \phi \cdot \mathbf{n} \rangle_{\partial \gamma}, \qquad \forall \phi \in\mathcal{D}(\Omega).
        \end{split}
    \end{align}

Now we apply Green's second identity to $-\Delta p$:   
    \begin{align*}
        -\langle \Delta \pD, \phi \rangle_{\Omega} &= -\langle \pD, \Delta \phi\rangle_{\Omega} \nolinebreak
                                                    + \langle \pD, \Delta \phi \cdot \mathbf{\nu} \rangle_{\partial\Omega} \nolinebreak
                                                    - \langle \nabla \pD \cdot \mathbf{\nu}, \phi \rangle_{\partial\Omega}\\
                                                    &= -\langle p, \Delta \phi\rangle_{\Omega} \nolinebreak \nolinebreak
                                                     + \langle p_D,  \nabla \phi \cdot \mathbf{\nu} \rangle_{\partial \Omega}.
    \end{align*}
By substituting \eqref{L_p_H2} into the above equation, we can get
    \begin{align}\label{L_p_D_1}
        -\langle \Delta \pD, \phi \rangle_{\Omega} = \langle f, \phi \rangle_{\Omega}  - \langle \pjump,  \nabla \phi \cdot \mathbf{n} \rangle_{\gamma}, \qquad \forall \phi \in \mathcal{D}(\Omega).
    \end{align}
Notice that
    \begin{align*}
        D_i\phi(\mathbf{y}) = D_i(\delta_0 * \phi(\mathbf{y})) = \phi *D_i\delta_0(\mathbf{y}) = \langle D_i\delta_0(\mathbf{y}-\mathbf{x}),\phi(\mathbf{x}) \rangle_\Omega,
    \end{align*}
we can convert $\langle \pjump,  \nabla \phi \cdot \mathbf{n} \rangle_{\gamma}$ into duality pairing over domain $\Omega$ by using the Fubini's theorem:
    \begin{align*}
        \langle \pjump,  \nabla \phi \cdot \mathbf{n} \rangle_{\gamma}  \nolinebreak
        &=\int_\gamma \pjump(\mathbf{y}) \nabla\phi(\mathbf{y}) \cdot \mathbf{n}(\mathbf{y}) \,d\sigma\\
        &= \int_\gamma \pjump(\mathbf{y})  \langle\nabla \delta_0(\mathbf{y}-\mathbf{x}) \cdot \mathbf{n}(\mathbf{y}), \phi(\mathbf{x})\rangle_\Omega  \,d \sigma\\
        & = \Big \langle \int_\gamma\pjump\nabla \delta_{\mathbf{x}} \cdot \mathbf{n} \,\ds, \phi(\mathbf{x}) \Big \rangle_\Omega.
    \end{align*}
    
Finally, we have
    \begin{align*}
        -\Big\langle \Delta \pD\vari, \phi\vari \Big\rangle_{\Omega} = \Big\langle \Big( f\vari - \int_\gamma\pjump\nabla \delta_{\mathbf{x}} \cdot \mathbf{n} \,\ds \Big),  \phi \vari \Big\rangle_{\Omega},  \qquad \forall \phi \in \mathcal{D}(\Omega),
    \end{align*}
and $-\Delta p$ can be expressed as
    \begin{align}\label{append}
        -\Delta \pD\vari = f\vari - \int_\gamma\pjump\nabla \delta_{\mathbf{x}} \cdot \mathbf{n} \,\ds, \qquad \forall \mathbf{x} \in \Omega.
    \end{align}

\end{appendices}

\end{document}